%% file: LTD.tex
\documentclass{amsart}
\input{headart}

\title{Local Tate duality for syntomic cohomology in positive characteristic}

\author{Christopher Lazda}
       \address{Department of Mathematics\\ Harrison Building \\ Streatham Campus
 \\ University of Exeter \\ North Park Road \\ Exeter \\ EX4 4QF \\ United Kingdom }
       \email{c.d.lazda@exeter.ac.uk}
       
\begin{document}

\begin{abstract}
Let $k$ be a local field of characteristic $p>0$.
We show that syntomic cohomology over $k$, with coefficients, and defined using Dodd's theory of $F$-gauges, satisfies local Tate duality.
We also prove that syntomic cohomology over the ring of integers of $k$, again with coefficients, is Lagrangian with respect to the duality pairing.
Finally, we show that our results here recover Shatz's duality and Milne's orthogonality results for the cohomology of finite flat group schemes of $p$-power order over $k$.
\end{abstract}

\maketitle 

\tableofcontents

\section*{Introduction}

Let $k$ be a non-archimedean local field, that is, a finite extension of either $\Q_p$ or $\F_p\lser{t}$ for some prime $p$.
Choose a separable closure $k^{\rm s}$ of $k$, and set $G_k={\rm Gal}(k^{\rm s}/k)$.
A fundamental result in the cohomology of such fields is local Tate duality \cite[Chapter I, Corollary 2.3]{Mil86}, which states that if $M$ is a finite (discrete) $G_k$-module, of order invertible in $k$, and we set
\[ M':= {\rm Hom}(M,\mu_\infty) \]
where $\mu_\infty$ is the Galois module of roots of unity in $k^{\rm s}$, then the trace map
\[ \rH^2(G_k,\mu_\infty)\to \Q/\Z \]
induces an isomorphism
\[ \rH^{2-i}(G_k,M') \isomto {\rm Hom}(\rH^{i}(G_k,M),\Q/\Z) \]
for all $i$.
If $\ell$ is a prime invertible in $k$, this implies a similar result for continuous $G_k$-representations in finitely generated $\Z_\ell$-modules.
Indeed, given such an $M$, we set $M^\vee:= {\bf R}{\rm Hom}_{\Z_\ell}(M,\Z_\ell)$ as a complex of continuous $\Z_\ell[G_k]$-modules, and define $M^\vee(1):=M^\vee\otimes_{\Z_\ell} \Z_\ell(1)$ where $\Z_\ell(1)=\lim{n} \mu_{\ell^n}(k^{\rm s})$ as usual.
Then the trace map
\[ {\bf R}\Gamma(G_k,\Z_\ell(1)[2]) \to \Z_\ell \]
induces an equivalence
\[ {\bf R}\Gamma(G_k,M^\vee(1)[2]) \isomto {\bf R}{\rm Hom}_{\Z_\ell}({\bf R}\Gamma(G_k,M),\Z_\ell) \]
in $\bD^b(\Z_\ell)$.

If ${\rm char}\,k=p$, and $M$ is a $G_k$-module of $p$-power order, then the na\"{\i}ve analogue of local Tate duality is false for trivial reasons.
Indeed, in this case there is no $p$-torsion in the $G_k$-module $\mu_\infty$, and so the above definition of $M'$ collapses to zero.
There is, of course, a natural way to fix this, by instead taking the dual in the category of finite flat group schemes.
That is, we consider $M$ instead as a finite \'etale group scheme over $k$, and then define
\[
M^D:=\underline{\rm Hom}_k(M,\G_m)
\]
to be its Cartier dual.
More generally, if $M$ is a (commutative) finite flat group scheme over $k$, we can define its Cartier dual $M^D$ in the same way.
The invariant map
\[
\rH^2(k,\G_m)\to \Q/\Z
\]
of local class field theory then induces a pairing 
\[
\rH^i(k,M)\times \rH^{2-i}(k,M^D)\to \Q/\Z
\]
This won't be a perfect pairing in general, because the cohomology groups involved need not be finite.
For example, $\rH^1(k,\Z/p)$, $\rH^1(k,\mu_p)$ and $\rH^1(k,\alpha_p)$ are all infinite.
However, these groups do have a natural locally compact topology (arising from the discrete valuation topology of $k$), and the main result of \cite{Sha64} is that the above pairing is a perfect pairing of locally compact abelian groups.

In the case of torsion prime to ${\rm char}\,k$, one of the main sources of Galois representations is the \'etale cohomology of varieties over $k$.
It has long been understood, however, that finite flat group schemes can only account for a small part of the $p$-adic cohomology of varieties over a field of positive characteristic (roughly speaking, the part with slopes in the interval $[0,1]$).
Hence one would like a more general version of local Tate duality over a positive characteristic local field, which would apply to a category of coefficients broad enough to accommodate the full $p$-adic cohomology of varieties over $k$.

One framework for these sort of duality results is the relatively perfect site introduced by Kato in \cite{Kat86}. Indeed, given a smooth variety $f\from X\to{\rm Spec}\,k$, Kato defines complexes ${\bf R}f_*\nu_n(r)$ on the relatively perfect site of $k$, whose cohomology can be shown to compute the absolute syntomic cohomology groups ${\bf R}\Gamma_{\rm Syn}(X,\Z/p^n(r))$.\footnote{Actually, Kato works over a perfect field, but his results should extend to the case of a field with a finite $p$-basis.} Kato then proves a duality result for ${\bf R}f_*\nu_n(r)$ as a sheaf on the relatively perfect site, and results of Suzuki \cite{Suz22,Suz24} can then be used to deduce a topological version of local Tate duality, as applied to the full relative relative cohomology complex of $X/k$.

Our aim in this article is to give an alternative approach to such duality results, working instead in with Dodd's theory of $F$-gauges \cite{Dod22}, a category large enough to accomodate all finite flat group schemes over $k$, as well as the relative cohomology of all smooth and proper varieties over $k$.
These $F$-gauges are graded modules over an analogue of Berthelot's ring of differential operators, and should be seen as a concrete realisation of the of prismatic $F$-gauges introduces by Bhatt--Lurie and Drinfeld \cite{Bha22,BL22,Dri20} (inspired by earlier work of Fontaine--Jannsen \cite{FJ21}).
In particular, we expect that the category of ``lisse $F$-gauges'' that we define should be equivalent to the category of perfect complexes on the syntomification $k^{\rm Syn}$ of $k$, although we do not attempt to prove this comparison here.
The fact that this category is a reasonable setting for duality results in $p$-adic cohomology is indicated by the recovery of local Tate duality in mixed characteristic via prismatic methods in \cite{Bha22}, as well as the geometric duality results for prismatic cohomology proved in \cite{Tan24}.

As is usual when studying the theory of arithmetic $\mathcal{D}$-modules, the construction depends on choosing a Cohen ring $C$ for $k$, and a Frobenius lift $\sigma$ on $C$, although the resulting category will be independent of these choices.
Using Dodd's category, or rather, a Frobenius descent of Dodd's category, we can define a relatively concrete version of syntomic cohomology ${\bf R}\Gamma_{\rm Syn}(k,\,\cdot\,)$ by taking the mapping complex from the unit object $C^\star$.
Similarly, any $F$-gauge $M$ satisfying suitable finiteness conditions has a dual $M^\vee$, given by taking the internal hom to $C^\star$ in the category of $F$-gauges.
The Tate twist $C^\star(1)$ makes sense in Dodd's category, and we can explicitly construct a trace map
\[
{\bf R}\Gamma_{\rm Syn}(k,C^\star(1)[2]) \to \Z_p
\]
which is essentially given by the familiar residue map on differentials.
As with finite flat group schemes, these syntomic cohomology groups are not in general finitely generated over $\Z_p$, and so any form of Tate duality will need to take account of the natural topology on syntomic cohomology. To allow for cleaner duality results, we replace these topologies with condensed structures in the sense of \cite{Sch26,CS23}.
For this $\mathcal{D}$-module theoretic syntomic cohomology, viewed as taking values in condensed $\Z_p$-modules, we can then prove the following form of local Tate duality.

\begin{theoremu}[\ref{theo: main1}]
Let $k$ be a local field of characteristic $p$, and $M$ a lisse $F$-gauge over $k$, with dual $M^\vee$. Then the trace map induces an equivalence
\[  
{\bf R}\Gamma_{\rm Syn}(k,M^\vee(1)[2]) \isomto  {\bf R}{\rm Hom}_{\Z_p}({\bf R}\Gamma_{\rm Syn}(k,M),\Z_p) 
\]
in the derived category of condensed $\Z_p$-modules.
\end{theoremu}

The condition of being lisse is a natural finiteness condition on $F$-gauges, essentially analogous to a complex of continuous $G_k$-representations in $\Z_\ell$-modules being bounded with finitely generated cohomology.
We fully expect that $F$-gauges associated to smooth and proper varieties over $k$ will satisfy this condition (although we do not address this question here).
As it turns out, the main difficulty in proving Theorem \ref{theo: main1} lies in setting up the appropriate machinery, once this has been done the proof itself is a pleasingly direct reduction to Pontryagin duality for finite dimensional $k$-vector spaces, via an explicit description of syntomic cohomology complexes.

We should emphasise that our result here is a purely $\mathcal{D}$-module theoretic one in the framework of \cite{Dod22}; we do not claim anything about the coherent cohomology of the stack $k^{\rm Syn}$.
However, the expected comparison results between these two approaches should eventually translate Theorem \ref{theo: main1} into a result on syntomic cohomology in the sense of \cite{Bha22}.

We leave the construction of $F$-gauges of `geometric origin' for future work; the only natural sources of $F$-gauges we consider here are virtual $F$-crystals, finite flat group schemes of $p$-power order, and $p$-divisible groups.
In each case, we prove that syntomic cohomology computes the expected cohomology groups: for virtual $F$-crystals this is homotopy fixed points of the divided Frobenius (see Theorem \ref{theo: syn cohom complex 2}), and for finite flat group schemes or $p$-divisible groups this is just flat cohomology.
In particular, we show that Theorem \ref{theo: main1} recovers Shatz's duality for finite flat group schemes of $p$-power order as a special case.

The fact that we do not attempt to construct gauges of geometric origin also means that we are unable here to prove precise comparisons with the duality theorems of Kato and Suzuki, although it is clear what the relationship should be: if $X$ is a smooth and proper scheme over $k$, then gauge-theoretic Poincar\'e duality for $X/k$, together with Theorem \ref{theo: main1} should essentially recover the absolute duality results for $X$ obtained by combining \cite{Kat86} with \cite[Proposition 4.7]{Suz24}. It is also worth pointing out that the framework of Kato and Suzuki could be used as an alternative to the condensed Pontryagin duality used here, and one which could lead to finer results on the resulting condensed structure on syntomic cohomology.

Finally, we show the expected Lagrangianity of syntomic cohomology over the ring of integers $k^\circ$ of $k$ with respect to local Tate duality, which should be seen as a function field analogue of \cite[Corollary 6.4.6]{Bha22}.
For this, we need to choose a $p$-adically complete subring $C^\circ \subset C$, stable under $\sigma$, satisfying $C^\circ/p=k^\circ$.
Given $(C,\sigma)$, it may not always be possible to choose such a $C^\circ$, but there are always choices of $(C,\sigma)$ for which such a $C^\circ$ exists.

\begin{theoremu}[\ref{theo: lagrange}]
For any lisse $F$-gauge
 $M^\circ$ over $k^\circ$, with base change $M$ to $k$, the sequence
 \begin{equation*}
 {\bf R}\Gamma_{\rm Syn}(k^\circ,M^\circ) \to  {\bf R}\Gamma_{\rm Syn}(k,M)\to  {\bf R}{\rm Hom}_{\Z_p}({\bf R}\Gamma_{\rm Syn}(k^\circ,M^{\circ,\vee}(1)[2]),\Z_p) \overset{+1}{\lto}
 \end{equation*}
induced by the duality pairing is a fibre sequence.
\end{theoremu}

Again, if $G$ is a finite flat group scheme over $k^\circ$, this recovers existing results \cite[Corollary III.7.2]{Mil86} on the orthogonality of $\rH^i(k^\circ,G)$ and $\rH^{2-i}(k^\circ,G^D)$ with respect to the local duality pairing. And for gauges arising from the cohomology of smooth and proper $k^\circ$-schemes, this should essentially recover a special case of the duality results proved in \cite{Zha18}.

The methods used to prove Theorems \ref{theo: main1} and \ref{theo: lagrange} generalise in a relatively straightforward manner to show that Poitou--Tate duality for syntomic cohomology over global function fields. This is the topic of work currently in preparation.

\subsection*{Overview}

In \S\ref{sec: gauge} we give the basic definitions concerning condensed gauges over $k$, as well as the various natural finiteness conditions on such objects, in particular the important condition of being lisse. We work at a general level $m\geq 0$, with level $m=1$ recovering the framework of \cite{Dod22}, and level $m=0$ being the most useful for cohomological calculations.

In \S\ref{sec: cond F desc} we explain how to deduce results on Frobenius descent for condensed $\mathcal{D}$-modules in a relatively straightforward manner from existing results in the literature on Frobenius descent for abstract $\mathcal{D}$-modules.
We use this to define $F$-gauges, and also prove a result on Frobenius descent for ($F$-)gauges.
This shows that the category of quasi-nilpotent ($F$-)gauges is independent of the level $m$.

In \S\ref{sec: symm} we define the symmetric monoidal structure on the category of ($F$-)gauges, and give an explicit formula for the dual of a(n $F$-)gauge.

In \S\ref{sec: indep} we show that various categories of gauges are independent of the choice of $(C,\sigma)$ up to canonical equivalence.

In \S\ref{sec: cohom} we define the various cohomology groups associated to an $F$-gauge of level $0$, namely Nygaard, de\thinspace Rham and syntomic, and give explicit complexes calculating these cohomology groups.

In \S\ref{sec: vfc etc} we show how to construct $F$-gauges from virtual $F$-crystals, $p$-divisible groups, and finite flat group schemes of $p$-power order over $k$.
We also show that syntomic cohomology computes the expected cohomology groups in each case.

In \S\ref{sec: PD TM} we formulate Pontryagin duality for finite dimensional $k$-vector spaces and finitely generated $C$-modules in the condensed framework, and construct the trace morphism in syntomic cohomology.

In \S\ref{sec: core} we prove our main result on local Tate duality, by combining the explicit complexes representing syntomic cohomology with the Pontryagin duality results from \S\ref{sec: PD TM}.

Finally, in \S\ref{sec: lagrange} we show that syntomic cohomology over $k^\circ$ is Lagrangian with respect to the duality pairing; again, this essentially reduces to Pontryagin duality for finitely generated $k^\circ$-modules.

\subsection*{Acknowledgements}
The author would like to thank Christopher Dodd, Oliver Gregory, Andreas Langer, Bernard Le Stum, Alexei Skorobogatov, Takashi Suzuki and Tianchen Zhao for helpful conversations regarding the contents of this article. Further acknowledgements to be added after the refereeing process.

\subsection*{AI Declaration}
Claude Opus 5 was used in the preparation of this article, mainly to perform literature searches, and for proof-reading earlier versions of the article. It made the following mathematically significant contributions:
\begin{enumerate}
\item the suggestion of the article \cite{Shi15} as a source for Frobenius descent to level $-1$;
\item the suggestion of the article \cite{TV24} as a source for the computation of flat cohomology via syntomic complexes;
\item the example of a non-gauge-like Dieudonn\'e module in Example \ref{exa: non gauge-like};
\item the example of a $(p,t)$-torsion $F$-gauge over $k^\circ$ in Remark \ref{rem: tors gauge}, suggested as a replacement for a previous example which did not admit an $F$-structure.
\end{enumerate}
All other mathematical content is the author's own.

 \subsection*{Notations and conventions}

\begin{itemize}
\item
Throughout, $k$ will denote a local field of characteristic $p$, with residue field $\F$ a finite extension of $\F_p$.
We will let $C$ denote a Cohen ring for $k$, and $\sigma$ a Frobenius lift on $C$ (the dependence on the choices of $C$ and $\sigma$ is discussed in \S\ref{sec: indep}). 
\item
We use cohomological indexing conventions for complexes, which will therefore be written as
\[ \ldots \overset{d}{\lto} M^{q-1} \overset{d}{\lto} M^{q} \overset{d}{\lto} M^{q+1} \overset{d}{\lto} \ldots \]
with the differential of degree $+1$.
Thus `bounded above' will mean that $M^q=0$ for $q\gg0$ \&c.
Shifts are then defined so that $M[n]^q=M^{q+n}$, in particular $\rH^n(M)=\rH^0(M[n])$. 
We say that a complex is static if it only has non-zero cohomology in degree $0$.
\item
We will use throughout the condensed sets introduced in \cite{Sch26} in place of topological spaces.
In fact, it will be enough for us to work in the light condensed framework described in \cite{CS23}.
Thus, we will consider the site of light profinite sets to be the category of countable limits of finite sets, equipped with the topology of finite collections of jointly surjective maps.
The bound on the cardinality of the indexing category of the limit ensures that the underlying category is essentially small, and a condensed gadget will then be a sheaf of gadgets on the site of light profinite sets; note that we call these objects condensed gadgets and not light condensed gadgets.
\item All finite flat group schemes will be assumed commutative.
\end{itemize}

\section{The category of gauges}
\label{sec: gauge}

Since $k$ is a ring in characteristic $p>0$, its syntomification $k^{\rm Syn}$ is defined in \cite{Bha22} as a stack over $\spf{\Z_p}$.
Syntomic cohomology is then given by the cohomology of the stack $k^{\rm Syn}$, with coefficients in vector bundles (or more generally, perfect complexes).
A more down-to-earth description of the syntomic cohomology of $k$ is conjecturally given by the theory of $F$-gauges introduced in \cite{Dod22}, and in this section we recall the key definitions concerning gauges over $k$, Frobenius structures on gauges will then be discussed in \S\ref{sec: cond F desc} below.

Perhaps the main difference between our approach and that of Dodd is that (eventually) we will want to work `one level down', in order to simplify explicit cohomological calculations.
Thus we develop the theory of gauges of arbitrary level $m$, even though we are ultimately mostly interested in the case $m=0$. The case $m=1$ is needed in order to clarify the link with \cite{Dod22}, and it is no harder to treat the case of general $m$ than it is to give a unified treatment of the cases $m=0$ and $m=1$.
With the minor exception of \S\ref{sec: vfc etc}, the reader is welcome to take $m=0$ throughout, which we shall do ourselves from \S\ref{sec: PD TM} onwards.

As well as dealing with gauges of different levels, we will also need a topological (or rather, condensed) enhancement of the approach used in \cite{Dod22}, in order to handle non-finitely generated syntomic cohomology groups.
Thus, rather than considering abstract gauges over $k$ as in \cite{Dod22}, we will instead consider what might be called `light condensed gauges'.
Since we will only ever need to use light condensed sets, we will generally drop the adjective `light'.

To begin with, then, we pick a Cohen ring $C$ for $k$, together with a lift $\sigma\from C\to C$ of the absolute Frobenius of $k$.
We let $\F$ denote the (finite) residue field of $k$, and $W=W(\F)$ the ring of Witt vectors of $\F$.
Thus $C$ is canonically a $W$-algebra.
Since any uniformiser of $k$ is a $p$-basis, the module $\Omega^1_{k}=\Omega^1_{k/\F}$ of differentials of $k$ is free of rank $1$, and hence the module $\widehat{\Omega}^1_C=\widehat{\Omega}^1_{C/W}$ of $p$-adically continuous differentials of $C$ is also free of rank $1$.
In particular, all the usual descriptions of Berthelot's rings of arithmetic differential operators apply, even though $k$ is not of finite type over $\F$, see for example \cite{CV22}. 

\subsection{The ring $\bm{\widehat{\mathcal{D}}^{(m-1,m),\star}_{C}}$}

Each ring $C/p^n$ has a natural locally compact topology: any uniformiser lift $t\in C$ induces an isomorphism
\[
C/p^n\cong W/p^n\lser{t}
\]
and the resulting $t$-adic topology on $C/p^n$ is independent of the choice of $t$. 
We can convert this into a condensed structure in the usual way, by setting
\[
C/p^n(S) = {\rm Map}_{\rm cts} (S,C/p^n)
\]
for any light profinite set $S$.
We then set $C=\lim{n}C/p^n$ as a condensed ring.
It is straightforward to check that $\rH^i(S,C)=0$ for any light profinite set $S$ and any $i>0$, and that the sequence
\[
0 \to C(S) \overset{p^n}{\lto} C(S)\lto C/p^n(S)\to 0
\]
is exact for all $n\geq 1$.
We view the various rings $W$, $\Z_p$, $W/p^n$ and $\Z/p^n$ as condensed rings in an analogous way.

For any $m,r\in \N$, the finite free $C/p^n$-module $\mathcal{D}^{(m)}_{C/p^n,r}$ of differential operators of order $\leq r$ and level $m$ is equipped with a natural $t$-adic topology, which we convert into a condensed structure in the same manner.
Thus we obtain the non-commutative condensed rings
\[
\mathcal{D}^{(m)}_{C/p^n}:=\colim{r} \mathcal{D}^{(m)}_{C/p^n,r}\qquad\widehat{\mathcal{D}}^{(m)}_{C}:=\lim{n} \mathcal{D}^{(m)}_{C/p^n}  \]
For any choice of $t$ as above, we have the section-wise analogue
\[
\widehat{\mathcal{D}}^{(m)}_{C}(S) = \left\{\left. \sum_{j=0}^\infty a_j\partial^{\langle j \rangle_{(m)}} \,\right\vert\, a_j\in C(S),\;a_j\to 0\; \text{uniformly }p\text{-adically}\right\}
\]
of the usual concrete description of elements of $\widehat{\mathcal{D}}^{(m)}$, where $\partial=\frac{\d}{{\d}t}$, and the $m$-PD differential operators $\partial^{\langle j \rangle_{(m)}}$ are as in \cite[\S2]{Ber96a}.
The natural left action of $\widehat{\mathcal{D}}^{(m)}_{C}$ on $C$ extends canonically to the condensed versions of these objects.

Throughout the rest of this article, the notation $C$, $\widehat{\mathcal{D}}^{(m)}_C$, \&c.\ will refer to the condensed rings constructed above, even if we may occasionally abuse notation and use the same notation for the underlying ring.
If we need to be careful about the distinction, we will use the standard notation $C(*)$, $\widehat{\mathcal{D}}^{(m)}_C(*)$ \&c.\ for the underlying ring of a condensed ring. In particular, a $C$-module (resp.\ $\widehat{\mathcal{D}}_C^{(m)}$-module, \&c.) will be implicitly understood to be a condensed $C$-module, rather than an abstract $C(*)$-module (resp.\ $\widehat{\mathcal{D}}^{(m)}_C(*)$-module, \&c.). We will say that a condensed module over any of these rings is finitely generated if it is in the essential image of the pullback functor from finitely generated abstract  (that is, uncondensed) modules.

We will also need the ring of differential operators of level $-1$, as introduced in \cite{Shi15}. The following quick-and-dirty definition will be good enough for us.

\begin{definition}
We define $\widehat{\mathcal{D}}^{(-1)}_{C}$ to be the $p$-completion of the subring of $ \widehat{\mathcal{D}}^{(0)}_{C}$ generated by $C$ and $p\partial$.
For each $n\geq 1$, we define $\mathcal{D}^{(-1)}_{C/p^n}$ to be $(\widehat{\mathcal{D}}^{(-1)}_{C})/p^n$.
\end{definition}

One easily checks that this doesn't depend on the choice of uniformiser lift $t$. We extend the definition of the $m$-PD differential operators to the case $m=-1$ by setting $\partial^{{\langle j\rangle}_{(-1)}}:=(p\partial)^j$, we will also use the notation $\theta:=p\partial$, as in \cite{Shi15}.

\begin{definition}
We denote by $\bD(\widehat{\mathcal{D}}^{(m)}_C)$ the $\infty$-category of derived $p$-complete complexes of left $\widehat{\mathcal{D}}^{(m)}_C$-modules, and by $\bD(\mathcal{D}^{(m)}_k)$ the $\infty$-category of complexes of left $\mathcal{D}^{(m)}_k$-modules.
\end{definition}

We make this definition in order to highlight the fact that such objects are always assumed derived $p$-complete, even if this is not clear from the notation. We also emphasise that these are derived categories of \emph{condensed} objects, and that they are understood to be derived categories in the $\infty$-categorical sense.

\begin{definition}
Following \cite{Dod22}, for any $m\geq 0$, we define the $\widehat{\mathcal{D}}^{(m-1)}_C$-submodule
\[
\widehat{\mathcal{D}}^{(m-1,m),i}_C := \left\{\left. P\in \widehat{\mathcal{D}}^{(m)}_C \,\right\vert\, p^iP\in \widehat{\mathcal{D}}^{(m-1)}_C \right\}
\]
and the corresponding graded ring
\[
\widehat{\mathcal{D}}^{(m-1,m),\star}_C:=\bigoplus_i \widehat{\mathcal{D}}^{(m-1,m),i}_C
\]
\end{definition}

We will mostly be interested in the cases $m=0,1$, and eventually just in the case $m=0$. For each $n\geq 1$, we also set
\[
\mathcal{D}^{(m-1,m),\star}_{C/p^n}:=\widehat{\mathcal{D}}^{(m-1,m),\star}_C/p^n
\] 
to be the mod $p^n$ version of this graded ring. There are natural maps
\[
f\from\widehat{\mathcal{D}}^{(m-1,m),i}_C\to \widehat{\mathcal{D}}^{(m-1,m),i+1}_C
\]
given by the inclusion, and 
\[
v\from \widehat{\mathcal{D}}^{(m-1,m),i}_C\to \widehat{\mathcal{D}}^{(m-1,m),i-1}_C
\]
given by multiplication by $p$.
These can be viewed as multiplication by commuting elements $f\in \widehat{\mathcal{D}}^{(m-1,m),1}_C(*)$, $v\in \widehat{\mathcal{D}}^{(m-1,m),-1}_C(*)$ satisfying $fv=p$.

A concrete description of $\widehat{\mathcal{D}}^{(m-1,m),\star}_C$ can then be given exactly as in \cite[\S3.1]{Dod22}.
We fix a uniformiser lift $t\in C(*)$, with corresponding derivation $\partial$. If $P\in\widehat{\mathcal{D}}^{(m-1)}_C$, then we define $[\partial^{[p^m]},P]_{(m-1)}\in \widehat{\mathcal{D}}^{(m-1)}_C\subset \widehat{\mathcal{D}}^{(m)}_C$ to be the commutator of $\partial^{[p^m]}$ and $P$, that this does indeed lie inside $\widehat{\mathcal{D}}^{(m-1)}_C$ can be checked using the standard relations for $m$-PD differential operators \cite[\S2]{Ber96a}. 

We then consider the non-commutative ring $\mathcal{A}^{(m)}$ generated over $\widehat{\mathcal{D}}^{(m-1)}_C$ by the commuting variables $f,v$ and $\partial^{[p^m]}=\partial^{\langle p^m\rangle_{(m)}}$, whose commutation relations with local sections of $\widehat{\mathcal{D}}^{(m-1)}_C$ are given by the rules:
\begin{enumerate} 
\item
$[f,P]=[v,P]=0$ for all $P\in \widehat{\mathcal{D}}^{(m-1)}_C$;
\item
$[\partial^{[p^m]},P]=f[\partial^{[p^m]},P]_{(m-1)}$ for all $P\in \widehat{\mathcal{D}}^{(m-1)}_C$.
\end{enumerate}  
In particular, for any light profinite set $S$, every element of $\mathcal{A}^{(m)}(S)$ can be written uniquely as a finite sum of the form $\sum_{i,j,k}\alpha_{i,j,k}(\partial^{[p^m]})^if^jv^k$ for some $\alpha_{i,j,k}\in \widehat{\mathcal{D}}^{(m-1)}_C(S)$, and $\mathcal{A}^{(m)}$ can be graded by giving $f$ and $\partial^{[p^m]}$ degree $1$, and $v$ degree $-1$. 

\begin{lemma}
Let $c_m\in\Z_p^\times$ be $(p-1)!$ if $m\geq 1$, and $1$ if $m=0$. Then the left ideal $I$ of $\mathcal{A}^{(m)}$ generated by the degree $0$ global sections 
\[
fv-p,\;\partial^{\langle p^m \rangle_{(m-1)}}-c_mv\partial^{[p^m]}\in \mathcal{A}^{(m),0}(*)
\]
is a centralising ideal, and in particular a two-sided ideal.
The graded ring $\widehat{\mathcal{D}}^{(m-1,m),\star}_C$ is isomorphic to $\mathcal{A}^{(m)}/I$.
\end{lemma}

\begin{proof}
The element $fv-p$ is in the centre of $\mathcal{A}^{(m)}$, and the element $\partial^{\langle p^m \rangle_{(m-1)}}-c_mv\partial^{[p^m]}$ commutes with $f,v$ and $\partial^{[p^m]}$, so to prove that $I$ is a centralising ideal it is enough to show that if $P\in \widehat{\mathcal{D}}^{(m-1)}_C$, then $[\partial^{\langle p^m \rangle_{(m-1)}}-c_mv\partial^{[p^m]},P]\in I$.
But we can calculate
\begin{align*}
[\partial^{\langle p^m \rangle_{(m-1)}}-c_mv\partial^{[p^m]},P] 
&=  [\partial^{\langle p^m \rangle_{(m-1)}},P]-c_mv[\partial^{[p^m]},P] \\
&=pc_m[\partial^{[p^m]},P]_{(m-1)}-c_mfv[\partial^{[p^m]},P]_{(m-1)} \\
&=c_m[\partial^{[p^m]},P]_{(m-1)}(p-fv) \in I
\end{align*} 
That $\widehat{\mathcal{D}}^{(m-1,m),\star}_C\cong \mathcal{A}^{(m)}/I$ can be shown entirely similarly to the case $m=1$ handled in \cite[\S3.1]{Dod22}. 
\end{proof}

\begin{remark}
In \cite{Dod22} the ring $\widehat{\mathcal{D}}^{(0,1),\star}_C$ is described by the weaker relation $f\partial^p=p!\partial^{[p]}$.
One easily checks from the definition that the stronger relation $\partial^p=(p-1)!v\partial^{[p]}$ holds.
\end{remark}

%\begin{corollary} \label{cor: gl sec}
%For all $m\geq 0$, the graded ring $\widehat{\mathcal{D}}^{(m-1,m),\star}_C$ is left and right coherent.
%The functor of global sections is an equivalence of categories between coherent graded $\widehat{\mathcal{D}}^{(m-1,m),\star}_C$-modules and finitely generated graded $\widehat{\mathcal{D}}^{(m-1,m),\star}_C(*)$-modules.
%\end{corollary}

\subsection{Condensed gauges over $\bm{k}$}

We can now define our main objects of study.

\begin{definition}
\label{defn: gauges}
A (condensed) gauge of level $m$ over $k$ is a derived $p$-complete complex of graded left $\widehat{\mathcal{D}}^{(m-1,m),\star}_C$-modules.
The $\infty$-category of such objects will be denoted $\bD(\widehat{\mathcal{D}}^{(m-1,m),\star}_C)$.
A (condensed) $p$-torsion gauge of level $m$ over $k$ is a complex of graded left $\mathcal{D}^{(m-1,m),\star}_k$-modules, and the $\infty$-category of such objects will be denoted $\bD(\mathcal{D}^{(m-1,m),\star}_k)$.
\end{definition}

Here (derived) $p$-completeness should be understood in the graded sense, in other words as $p$-completeness of each graded piece.
We will generally drop the adjectives `derived' and `graded' when talking about $p$-completeness; this will essentially always be meant in the derived sense, and whether it is meant in the ordinary or graded sense should be clear from context.
We will also generally drop the adjective `condensed', and refer simply to gauges over $k$.
We will use all the usual adjectives (bounded above/below, perfect \&c.) to isolate particular subcategories of ${\bD}(\widehat{\mathcal{D}}^{(m-1,m),\star}_C)$.
%, for example, we will refer to the category of bounded, coherent gauges of level $m$, denoted $\bD^b_\coh(\widehat{\mathcal{D}}^{(m-1,m),\star}_C)$, or the category of perfect $p$-torsion gauges of level $m$, denoted $\bD_\perf(\mathcal{D}^{(m-1,m),\star}_k)$.
One might argue that what we have called gauges of level $m$ would more naturally be called gauges of level $m-1$, we have chosen our convention so that Nygaard cohomology (see \S\ref{sec: cohom} below) is computed in the expected way at level $0$, see Proposition \ref{prop: cohom explicit 1}.

\begin{remark}
\label{rem: cond abs}
Since $\widehat{\mathcal{D}}^{(m-1,m),\star}_C$ is derived $p$-complete, there is a t-exact equivalence of categories between $\bD_{\rm perf}(\widehat{\mathcal{D}}^{(m-1,m),\star}_C)$ and the $\infty$-category of perfect complexes of $\widehat{\mathcal{D}}^{(m-1,m),\star}_C(*)$-modules.
Thus, although Definition \ref{defn: gauges} is rather naive (for example, we do not attempt to define a more sophisticated analytic ring structure on $\widehat{\mathcal{D}}^{(m-1,m),\star}_C$), the elementary category of gauges defined here will be good enough for our purposes.
\end{remark}

Informally, we think of a gauge of level $m$ as a diagram of $p$-complete complexes of $\widehat{\mathcal{D}}^{(m-1)}_C$-modules
\[
\ldots  \overset{f}{\underset{v}{\rightleftarrows}} M^{i-1} \overset{f}{\underset{v}{\rightleftarrows}} M^i  \overset{f}{\underset{v}{\rightleftarrows}} M^{i+1} \overset{f}{\underset{v}{\rightleftarrows}} \ldots
\]
with $fv=vf=p$, together with $W$-linear maps $\partial^{[p^m]}\from M^i\to M^{i+1}$, satisfying $[\partial^{[p^m]},P]=f[\partial^{[p^m]},P]_{(m-1)}$ for all $P\in \widehat{\mathcal{D}}^{(m-1)}_C$, $[\partial^{[p^m]},f]=[\partial^{[p^m]},v]=0$, and $v\partial^{[p^m]}= \frac{1}{c_m}\partial^{\langle p^m \rangle_{(m-1)}}$. Then
\[
M^{-\infty} = \colim{v} M^i = M/(v-1)
\]
is a complex of left $\widehat{\mathcal{D}}^{(m-1)}_C$-modules, and 
\[
M^{+\infty} = \colim{f} M^i = M/(f-1)
\]
is a complex of left modules over the subring of $\widehat{\mathcal{D}}^{(m)}_C$ generated by $\widehat{\mathcal{D}}^{(m-1)}_C$ and $\partial^{[p^m]}$.
Thus the $p$-completion $\widehat{M}^{-\infty}$ is an object of $\bD(\widehat{\mathcal{D}}^{(m-1)}_C)$, and the $p$-completion $\widehat{M}^{+\infty}$ is an object of $\bD(\widehat{\mathcal{D}}^{(m)}_C)$. 
Note that $\partial^{[p^m]}$ does induce a map $M^{-\infty}\to M^{-\infty}$, but the relations in $\widehat{\mathcal{D}}^{(m-1,m),\star}_C$ show that this map is equal to $\frac{1}{c_m}\partial^{\langle p^m\rangle_{(m-1)}}$, so there is no induced $\widehat{\mathcal{D}}^{(m)}_C$-module structure on $\widehat{M}^{-\infty}$. We denote by $f^\infty\from M^i\to \widehat{M}^{+\infty}$ and $v^\infty\from M^i\to \widehat{M}^{-\infty}$ the canonical maps.

\begin{example}\leavevmode
\begin{enumerate}
\item
The most important example of a gauge of level $m$ is the unit gauge $C^\star:=C[f,v]/(fv-p)$, equipped with its natural $\widehat{\mathcal{D}}^{(m-1,m),\star}_C$-module structure coming from the standard left $\widehat{\mathcal{D}}^{(m)}_C$-module structure on $C$.
This gauge can be visualised by the diagram
\[
\ldots  \overset{p}{\underset{1}{\rightleftarrows}} C \overset{p}{\underset{1}{\rightleftarrows}} C  \overset{1}{\underset{p}{\rightleftarrows}} C \overset{1}{\underset{p}{\rightleftarrows}} \ldots
\]
with central term sitting in graded degree $0$.
\item
By reducing the unit gauge modulo $p$, we get the unit $p$-torsion gauge $k^\star:=k[f,v]/(fv)$, visualised by the diagram
\[
\ldots  \overset{0}{\underset{1}{\rightleftarrows}} k \overset{0}{\underset{1}{\rightleftarrows}} k  \overset{1}{\underset{0}{\rightleftarrows}} k \overset{1}{\underset{0}{\rightleftarrows}} \ldots
\]
\end{enumerate}
\end{example}

Note that any gauge (or $p$-torsion gauge) has an underlying complex of graded $C^\star$-modules, and any $p$-torsion gauge has an underlying complex of graded $k^\star$-modules. We can produce more examples of gauges using Tate twists. 

\begin{definition}
If $M$ is a gauge, and $j\in \Z$, we define the Tate twist $M(j)$ by shifting the grading according to the rule
\[
M(j)^i = M^{i-j}
\]
\end{definition}

\begin{example}
The Tate twist $C^\star(1)$ can be visualised by the diagram
\[
\ldots  \overset{p}{\underset{1}{\rightleftarrows}} C \overset{p}{\underset{1}{\rightleftarrows}} C  \overset{p}{\underset{1}{\rightleftarrows}} C \overset{1}{\underset{p}{\rightleftarrows}} \ldots 
\]
with central term sitting in graded degree $0$.
\end{example}

The functors 
\begin{align*}
\widehat{(\,\cdot\,)}^{+\infty}\from &\bD(\widehat{\mathcal{D}}^{(m-1,m),\star}_C)\to  \bD(\widehat{\mathcal{D}}^{(m)}_C)\\
\widehat{(\,\cdot\,)}^{-\infty}\from &\bD(\widehat{\mathcal{D}}^{(m-1,m),\star}_C)\to  \bD(\widehat{\mathcal{D}}^{(m-1)}_C)
\end{align*}	
both have right adjoints, denoted $j_+^{(m)}$ and $j_+^{(m-1)}$ respectively. Concretely, we have $(j_+^{(m)}M)^i=M$ for all $i$, with $f={\rm id}$, $v=p$, $\mathcal{\widehat{D}}_{C}^{(m-1)}$ acting on each $(j_+^{(m)}M)^i=M$ in the given way, and $\partial^{[p^m]}\from (j_+^{(m)}M)^i=M\to (j_+^{(m)}M)^{i+1}=M$ again being the given action. Similarly, we have $(j_+^{(m-1)}M)^i=M$ for all $i$, with $f=p$, $v=\mathrm{id}$, $\mathcal{\widehat{D}}_{C}^{(m-1)}$ acting on each $(j_+^{(m-1)}M)^i=M$ in the given way, and $\partial^{[p^m]}\from (j_+^{(m-1)}M)^i=M\to (j_+^{(m-1)}M)^{i+1}=M$ being the given action of $\frac{1}{c_m}\partial^{\langle p^m\rangle_{(m-1)}}$.

Note that the essential image of $j_+^{(m)}$ (resp.\ $j_+^{(m-1)}$) is given by the full subcategory of gauges $M$ for which each $f\from M^i\to M^{i+1}$ (resp.\ each $v\from M^i\to M^{i-1}$) is an equivalence. We also readily see that $(j^{(m)}_+M)^{+\infty}=M$, whilst $(j^{(m)}_+M)^{-\infty}=\colim{\times p}M=M[1/p]$. Hence after $p$-completing these become $M$ and $0$ respectively. Similarly, $(j_+^{(m-1)}M)^{+\infty}=M[1/p]$ and $(j^{(m-1)}_+M)^{-\infty}=M$, which again become $0$ and $M$ after $p$-completing. Be warned that the notation $j_+^{(m)}$ is potentially ambiguous, since it could be interpreted as taking values in either $\bD(\widehat{\mathcal{D}}^{(m-1,m),\star}_C)$ or $\bD(\widehat{\mathcal{D}}^{(m,m+1),\star}_C)$. It should always be clear from context which of the two is meant.

\subsection{Quasi-nilpotent and lisse gauges}

Quasi-nilpotence conditions for gauges can be expressed using the fact that, for all $m\geq 0$, the element $\partial^{\langle p^m\rangle_{(m-1)}}$ is central in both $\mathcal{D}^{(m-1)}_k$ and $\mathcal{D}^{(m-1,m),\star}_k$. 

\begin{definition} 
A $p$-torsion gauge of level $m$ over $k$, or a complex of left $\mathcal{D}_k^{(m-1)}$-modules, is called quasi-nilpotent if it is cohomologically $(\partial^{\langle p^m\rangle_{(m-1)}})^\infty$-torsion.
That is, each of its cohomology groups is the union of its $(\partial^{\langle p^m\rangle_{(m-1)}})^n$-torsion submodules.
A gauge of level $m$ over $k$, or an object of $\bD(\widehat{\mathcal{D}}_C^{(m-1)})$, is called quasi-nilpotent if its (derived) mod $p$ reduction is quasi-nilpotent.
In each case, we will denote by $\bD_{\rm qn}\subset \bD$ the full subcategory of the relevant $\infty$-category consisting of quasi-nilpotent objects. 
\end{definition}

In particular, we see that a gauge is quasi-nilpotent if and only if each of its graded pieces is so. 
Moreover, if $M$ is a $p$-torsion gauge, then $M/^{\bf L}p \simeq M\oplus M[1]$, so $M$ being quasi-nilpotent does not depend on whether we consider it as a gauge or as a $p$-torsion gauge.
Note also that if $M$ is quasi-nilpotent, then $\widehat{M}^{-\infty}$ is quasi-nilpotent as a complex of $\widehat{\mathcal{D}}_C^{(m-1)}$-modules, but $\widehat{M}^{+\infty}$ need not be quasi-nilpotent as a complex of $\widehat{\mathcal{D}}_C^{(m)}$-modules.

In general, the question of appropriate finiteness conditions on gauges is an extremely subtle one. Since we are working over a field $k$, we can avoid this problem and make the following simple definition.

\begin{definition}
A gauge of level $m$ is called lisse if the underlying complex of graded $C^\star$-modules is perfect.
We denote by $\bD_{\rm lisse}(\widehat{\mathcal{D}}^{(m-1,m),\star}_C)\subset \bD(\widehat{\mathcal{D}}^{(m-1,m),\star}_C)$ the full subcategory of lisse gauges.
\end{definition}

Note that any perfect complex of graded $C^\star$-modules is automatically $p$-complete, this is not a condition that needs to be checked separately.
We will also find the following slightly weaker condition useful.

\begin{definition}
A gauge of level $m$ is called quasi-lisse if it is bounded, and each graded piece is perfect as a complex of $C$-modules.
\end{definition}

We can show exactly as in \cite[Lemma 2.1.1]{FJ21} that a quasi-lisse gauge $M$ is lisse if and only if the maps $f\from M^i\to M^{i+1}$, $v\from M^{-i}\to M^{-(i+1)}$ are equivalences for large enough $i$.
For lisse gauges, $M^{\pm\infty}$ are both already perfect over $C$, and thus $p$-complete, however, for quasi-lisse gauges, $M^{\pm\infty}$ need not be perfect over $C$ or $p$-complete.
%It is also easy to check that the analogues of Corollary \ref{cor: gl sec} and Remark \ref{rem: cond abs} hold for (complexes of) condensed graded $C^\star$-modules.
%Thus, since the underlying ring of $C^\star$ is regular, any lisse gauge is perfect as a complex of graded $C^\star$-modules.

\begin{lemma}
\label{lemma: qn nilp}
Any lisse gauge is quasi-nilpotent.
\end{lemma}

\begin{proof} 
We reduce modulo $p$ and replace $M$ by its cohomology groups, which are necessarily finitely generated as $C^\star$-modules.
Then $f\from M^i\to M^{i+1}$ is an isomorphism for $i\gg 0$, and so the map $v\from M^{i+1}\to M^i$ is zero for $i\gg 0$. 

If we fix $i_0$, then for $n$ large enough and $x\in M^{i_0}$ we have
\[
\left(\partial^{\langle p^m\rangle_{(m-1)}}\right)^n(x) = (c_mv\partial^{[p^m]})^n(x) =v^n(c_m\partial^{[p^m]})^n(x) \in v^{n-1}v(M^{i_0+n}) = 0
\]
and so $\left(\partial^{\langle p^m\rangle_{(m-1)}}\right)^n=0$ on $M^{i_0}$.
Since this is true for all $i_0$, we deduce that $M$ is quasi-nilpotent.
\end{proof}

\section{Condensed Frobenius descent and \texorpdfstring{$F$}{F}-gauges}
\label{sec: cond F desc}

A key input into the theory of gauges is Berthelot's theory of Frobenius descent \cite{Ber00}, and we will need a condensed version of this. In non-negative levels, this can be achieved quickly and easily, although in level $-1$ we will need to work a little harder. 

Denote by $C'$ the base change of $C$ along the Frobenius $W\to W$, we will abuse notation and also write $\sigma$ for the induced $W$-linear map $C'\to C$. For each pair $n,m$ of non-negative integers, Berthelot constructs an abstract (that is, noncondensed) $(\mathcal{D}^{(m+1)}_{C/p^n},\mathcal{D}^{(m)}_{C'/p^n})$-bimodule $\sigma^{*}\mathcal{D}_{C'/p^n}^{(m)}$, as well as an abstract $(\mathcal{D}^{(m)}_{C'/p^n},\mathcal{D}^{(m+1)}_{C/p^n})$-bimodule $\sigma^{\flat}\mathcal{D}_{C'/p^n}^{(m)}$. The former is flat as an abstract right $\mathcal{D}^{(m)}_{C'/p^n}$-module, and the latter is flat as an abstract right $\mathcal{D}^{(m+1)}_{C/p^n}$-module. Moreover, there is a canonical isomorphism
\[
\sigma^{\flat}\mathcal{D}_{C'/p^n}^{(m)} \otimes_{\mathcal{D}^{(m+1)}_{C/p^n}} \sigma^{*}\mathcal{D}_{C'/p^n}^{(m)}\isomto \mathcal{D}_{C'/p^n}^{(m)}
\]
of abstract $\mathcal{D}^{(m)}_{C'/p^n}$-bimodules, as well as a canonical isomorphism
\[
\mathcal{D}^{(m+1)}_{C/p^n} \isomto  \sigma^{*}\mathcal{D}_{C'/p^n}^{(m)} \otimes_{\mathcal{D}^{(m)}_{C'/p^n}} \sigma^{\flat}\mathcal{D}_{C'/p^n}^{(m)}
\] 
of abstract $\mathcal{D}^{(m+1)}_{C/p^n}$-bimodules. Taking the limit in $n$ we obtain a $\widehat{\mathcal{D}}^{(m)}_{C'}$-flat abstract $(\widehat{\mathcal{D}}^{(m+1)}_C,\widehat{\mathcal{D}}^{(m)}_{C'})$-bimodule $\sigma^{*}\widehat{\mathcal{D}}_{C'}^{(m)}$, and a $\widehat{\mathcal{D}}^{(m+1)}_{C}$-flat abstract $(\widehat{\mathcal{D}}^{(m)}_{C'},\widehat{\mathcal{D}}^{(m+1)}_{C})$-bimodule $\sigma^{\flat}\widehat{\mathcal{D}}_{C'}^{(m)}$, together with equivalences
\begin{align}
\label{eqn: frob desc}
\sigma^{\flat}\widehat{\mathcal{D}}_{C'}^{(m)} \widehat{\otimes}^{\bf L}_{\widehat{\mathcal{D}}^{(m+1)}_{C}} \sigma^{*}\widehat{\mathcal{D}}_{C'}^{(m)} & \isomto \widehat{\mathcal{D}}_{C'}^{(m)}  \\
\widehat{\mathcal{D}}_{C}^{(m+1)} & \isomto \sigma^{*}\widehat{\mathcal{D}}_{C'}^{(m)} \widehat{\otimes}^{\bf L}_{\widehat{\mathcal{D}}^{(m)}_{C'}} \sigma^{\flat}\widehat{\mathcal{D}}_{C'}^{(m)} \nonumber
\end{align}
Note that in each of these tensor products, the relevant modules are finite flat, so the abstract tensor product agrees with the $p$-completed tensor product, and the derived tensor product agrees with the underived one. As groups, $\sigma^{*}\widehat{\mathcal{D}}_{C'}^{(m)}$ and $\sigma^{\flat}\widehat{\mathcal{D}}_{C'}^{(m)}$ are just finite direct sums of copies of $\widehat{\mathcal{D}}_{C'}^{(m)}$, and so we may endow them with a canonical condensed structure in such a way that they become condensed bimodules, and the isomorphisms \eqref{eqn: frob desc} become isomorphisms of condensed bimodules. Thus we obtain:

\begin{theorem}
\label{theo: cdF 1}
Suppose $m\geq 0$. The functor $\sigma^*=\sigma^{*}\widehat{\mathcal{D}}_{C'}^{(m)} \widehat{\otimes}^{\bf L}_{\widehat{\mathcal{D}}_{C'}^{(m)}}(\,\cdot\,)$ is an equivalence of categories
\[
\bD(\widehat{\mathcal{D}}^{(m)}_{C'}) \isomto \bD(\widehat{\mathcal{D}}^{(m+1)}_C)
\]
with quasi-inverse given by $\sigma^{\flat}\widehat{\mathcal{D}}_{C'}^{(m)} \widehat{\otimes}^{\bf L}_{\widehat{\mathcal{D}}_{C}^{(m+1)}}(\,\cdot\,)$.
\end{theorem}

\begin{remark}
Again, since the modules $\sigma^{*}\widehat{\mathcal{D}}_{C'}^{(m)}$ and $\sigma^{\flat}\widehat{\mathcal{D}}_{C'}^{(m)}$ are of finite type, we can use the abstract tensor product of condensed modules; we do not need to worry about $p$-adic completion.
\end{remark}

In level $-1$ the theory of Frobenius descent is slightly more complicated. As in \cite{Shi15}, we can still define the left and right pullback functors $\sigma^*$ and $\sigma^\flat$, however, the maps \eqref{eqn: frob desc} are no longer isomorphisms, and so these functors need not be equivalences. However, Shiho showed in \cite[Theorem 3.1]{Shi15} that they are equivalences provided we restrict to (non-condensed) quasi-nilpotent objects, which will enable us to deduce a similar statement in the condensed setting. 

We do so by setting, for each $m\geq -1$,
\[
\breve{\mathcal{D}}_{C'}^{(m)}:= \lim{n} \widehat{\mathcal{D}}_{C'}^{(m)}/(p,\partial^{\langle p^{m+1} \rangle_{(m)}})^n
\]
to be the $(p,\partial^{\langle p^{m+1} \rangle_{(m)}})$-completion of $\widehat{\mathcal{D}}_{C'}^{(m)}$, and
\[
\breve{\mathcal{D}}_{C}^{(m+1)}:= \lim{n} \widehat{\mathcal{D}}_{C}^{(m+1)}/(p,\partial^{\langle p^{m+2} \rangle_{(m+1)}})^n
\]
to be the $(p,\partial^{\langle p^{m+2} \rangle_{(m+1)}})$-completion of $\widehat{\mathcal{D}}_{C}^{(m+1)}$.
Thus any quasi-nilpotent complex of $\widehat{\mathcal{D}}_{C'}^{(m)}$-modules (resp.\ $\widehat{\mathcal{D}}_{C}^{(m+1)}$-modules) canonically extends to a complex over $\breve{\mathcal{D}}_{C'}^{(m)}$ (resp.\ $\breve{\mathcal{D}}_{C}^{(m+1)}$). 

Playing the same game as before, we get a $(\breve{\mathcal{D}}_{C}^{(m+1)},\breve{\mathcal{D}}_{C'}^{(m)})$-bimodule structure on $\sigma^*\breve{\mathcal{D}}_{C'}^{(m)}$, and a $(\breve{\mathcal{D}}_{C'}^{(m)},\breve{\mathcal{D}}_{C}^{(m+1)})$-bimodule structure on $\sigma^\flat\breve{\mathcal{D}}_{C'}^{(m)}:={\rm Hom}_{C'}(C,\breve{\mathcal{D}}_{C'}^{(m)})$.
The crucial computation is then the following.

\begin{proposition}
\label{prop: bimod isom}
Suppose $m\geq -1$. Then there is a canonical isomorphism
\[
 \breve{\mathcal{D}}^{(m+1)}_{C}\isomto\sigma^{*}\breve{\mathcal{D}}_{C'}^{(m)} \otimes_{\breve{\mathcal{D}}^{(m)}_{C'}} \sigma^{\flat}\breve{\mathcal{D}}_{C'}^{(m)} 
\]
of abstract $\breve{\mathcal{D}}^{(m+1)}_{C}$-bimodules.
\end{proposition}

\begin{proof} 
If $m\geq 0$, this follows from suitably completing \eqref{eqn: frob desc}.
If $m=-1$, this can be proved by working modulo $p^n$ for all $n$, where we proceed more or less exactly as in \cite[Proposition 2.5.2]{Ber00}, using Shiho's calculations of the PD envelopes of level $-1$ from \cite[\S2.1]{Shi15} as the key input.
Of course, strictly speaking Shiho's results do not apply in our setting, since $k$ is not smooth over its residue field $\F$.
 However, the fact that any uniformiser $t$ is a $p$-basis for $k$ means that the relevant calculations of divided power envelopes still work in our setting.

The crucial difference we encounter when running through the proof of \cite[Proposition 2.5.2]{Ber00} is that we must work throughout with the full PD envelopes $\mathcal{P}_{C/p^n,(0)}$ and $\mathcal{P}_{C'/p^n,(-1)}$, and we never pass to the quotients by the natural PD filtration.
This gives the required answer because
\[
\breve{\mathcal{D}}_{C/p^n}^{(0)}:=\breve{\mathcal{D}}_{C}^{(0)}/p^n \isomto {\rm Hom}_{C/p^n}(\mathcal{P}_{C/p^n,(0)},C/p^n)
\]
and similarly at level $-1$. Note that we do \emph{not} get the analogous isomorphism
\[
 \mathcal{D}^{(0)}_{C/p^n} \cong \sigma^{*}\mathcal{D}_{C'/p^n}^{(-1)} \otimes_{\mathcal{D}^{(-1)}_{C'/p^n}} \sigma^{\flat}\mathcal{D}_{C'/p^n}^{(-1)} 
\]
for finite order differential operators, because the analogue of \cite[Lemme 2.3.4]{Ber00} does \emph{not} hold at level $-1$. This means that the isomorphism
\[
C/p^n\otimes_{C'/p^n} \mathcal{P}_{C'/p^n,(-1)} \otimes_{C'/p^n} C/p^n\isomto \mathcal{P}_{C/p^n,(0)}
\]
does \emph{not} descend to an isomorphism of pro-objects
\[
\{ C/p^n\otimes_{C'/p^n} \mathcal{P}^r_{C'/p^n,(-1)} \otimes_{C'/p^n}C/p^n\}_r \cong \{\mathcal{P}^r_{C/p^n,(0)}\}_r \qedhere
\]
\end{proof}

Given the proposition, we deduce that $\sigma^{*}\breve{\mathcal{D}}_{C'}^{(m)}$ (resp.\ $\sigma^{\flat}\breve{\mathcal{D}}_{C'}^{(m)}$) is a finite projective left (resp.\ right) $\breve{\mathcal{D}}_{C}^{(m+1)}$-module, and that 
\[ 
\sigma^{\flat}\breve{\mathcal{D}}_{C'}^{(m)} \otimes_{\breve{\mathcal{D}}^{(m+1)}_{C}} \sigma^{*}\breve{\mathcal{D}}_{C'}^{(m)}\isomto \breve{\mathcal{D}}_{C'}^{(m)}
\]
is an isomorphism of $\breve{\mathcal{D}}^{(m)}_{C'}$-bimodules, exactly as in \cite{Ber00}. Exactly as before, we can promote both this and Proposition \ref{prop: bimod isom} to isomorphisms of \emph{condensed} bimodules. We therefore deduce:

\begin{theorem}
\label{theo: cdF 2}
Suppose $m\geq -1$. Then the functor
\[
\sigma^*=\sigma^{*}\breve{\mathcal{D}}_{C'}^{(m)} \widehat{\otimes}^{\bf L}_{\breve{\mathcal{D}}_{C'}^{(m)}}(\,\cdot\,)
\]
is an equivalence of categories
\[
\bD_{\rm qn}(\widehat{\mathcal{D}}^{(m)}_{C'}) \isomto \bD_{\rm qn}(\widehat{\mathcal{D}}^{(m+1)}_C)
\]
with quasi-inverse given by $\sigma^{\flat}\breve{\mathcal{D}}_{C'}^{(m)} \widehat{\otimes}^{\bf L}_{\breve{\mathcal{D}}_{C}^{(m+1)}}(\,\cdot\,)$.
\end{theorem}

\begin{remark}
As before, the tensor products here can be abstract or $p$-adic, derived or underived, the resulting functors are the same.
\end{remark}

Composing with the base change $C\to C'$ along the Frobenius on $W$ therefore gives a semilinear pullback functor
\[
\sigma^* \from \bD(\widehat{\mathcal{D}}^{(m)}_C) \to \bD(\widehat{\mathcal{D}}^{(m+1)}_C)
\]
which is an equivalence for all $m\geq 0$, and which induces an equivalence
\[ 
\sigma^* \from \bD_{\rm qn}(\widehat{\mathcal{D}}^{(m)}_C) \isomto \bD_{\rm qn}(\widehat{\mathcal{D}}^{(m+1)}_C) 
\]
for all $m\geq -1$. We can now define Frobenius structures on gauges.

\begin{definition} 
An $F$-structure (or Frobenius structure) on a gauge $M$ of level $m$ is an equivalence
\[
\phi\from \sigma^*\widehat{M}^{-\infty}\isomto \widehat{M}^{+\infty}
\] 
in $\bD(\widehat{\mathcal{D}}^{(m)}_C)$.
An $F$-gauge of level $m$ is then a gauge of level $m$ equipped with an $F$-structure, and the category $\bD_{F}(\widehat{\mathcal{D}}^{(m-1,m),\star}_C)$ of $F$-gauges of level $m$ is defined to be the (strong) equaliser of the diagram
\[
\bD(\widehat{\mathcal{D}}^{(m-1,m),\star}_C) \xrightrightarrows[\sigma^*\widehat{(\,\cdot\,)}^{-\infty}]{\widehat{(\,\cdot\,)}^{+\infty}} \bD(\widehat{\mathcal{D}}^{(m)}_C) 
\] 
\end{definition}

Adjectives such as bounded, lisse, quasi-nilpotent, \&c.\ are understood to apply to the underlying gauge. In particular, if $M$ is a quasi-nilpotent $F$-gauge, then $\widehat{M}^{-\infty}$ is quasi-nilpotent as a $\widehat{\mathcal{D}}^{(m-1)}_C$-module, and so $\widehat{M}^{+\infty}$ is quasi-nilpotent as a $\widehat{\mathcal{D}}^{(m)}_C$-module.

\begin{remark}
\label{rem: F vs Finv and independence of sigma}
In \cite{Dod22}, such objects are called $F^{-1}$-gauges, but since the notion of $F$-gauges used there will not be relevant for us, we decided to use the slightly simpler terminology here.
It is also worth pointing out that our notion of $F$-gauge is marginally more general than that considered in \cite{Dod22} (at least, for graded $p$-complete objects), since we require an equivalence $\sigma^*\widehat{M}^{-\infty}\isomto \widehat{M}^{+\infty}$ rather than an equivalence $\sigma^*M^{-\infty}\isomto \widehat{M}^{+\infty}$.
In other words, we do not require $M^{-\infty}$ to be derived $p$-complete.
For all objects of interest to us, this distinction will be irrelevant. 
\end{remark}

We can use essentially the same methods as above to prove a Frobenius descent result for ($F$-)gauges. Indeed, if $m\geq 0$, then we can use finite flatness of $\sigma$ to describe
\[
\sigma^*\widehat{\mathcal{D}}_C^{(m-1,m),i} = \left\{\left. P \in \sigma^*\widehat{\mathcal{D}}^{(m)}_C \,\right\vert\,  p^iP\in \sigma^*\widehat{\mathcal{D}}^{(m-1)}_C \right\}
\]
from which we deduce the existence of a canonical graded $(\widehat{\mathcal{D}}_C^{(m,m+1),\star},\widehat{\mathcal{D}}_C^{(m-1,m),\star})$-bimodule structure on $\sigma^*\widehat{\mathcal{D}}_C^{(m-1,m),\star}$.
This results in a functor $\sigma^*$ from graded left $\widehat{\mathcal{D}}_C^{(m-1,m),\star}$-modules to graded left $\widehat{\mathcal{D}}_C^{(m,m+1),\star}$-modules, which derives trivially, and which on complexes preserves both $p$-completeness and quasi-nilpotence.
Moreover, on each graded piece $M^i$, or after passing to $\widehat{M}^{\pm\infty}$, this recovers the standard Frobenius base change functors. We therefore obtain functors
\begin{equation}
\label{eqn: F bc gauge}
\sigma^* \from\bD_{(F)}(\widehat{\mathcal{D}}^{(m-1,m),\star}_C)\to \bD_{(F)}(\widehat{\mathcal{D}}^{(m,m+1),\star}_C) 
\end{equation}
from ($F$-)gauges of level $m$ to those of level $m+1$. 
 
\begin{theorem}
\label{theo: F desc gauge}
If $m\geq 1$, then the functor \eqref{eqn: F bc gauge} is an equivalence of categories. If $m=0$, then \eqref{eqn: F bc gauge} induces an equivalence of categories
\[
\bD_{(F),{\rm qn}}(\widehat{\mathcal{D}}^{(-1,0),\star}_C)\isomto \bD_{(F),{\rm qn}}(\widehat{\mathcal{D}}^{(0,1),\star}_C) 
\]
\end{theorem}

\begin{proof}
It suffices to prove the claim for the categories of objects without Frobenius structure, and in fact it suffices to prove the analogous claim for $p$-torsion gauges, that is, complexes of graded left modules over $\mathcal{D}^{(m-1,m),\star}_k$.
Indeed, a straightforward $p$-completeness argument applied to ${\bf R}{\rm Hom}$'s then shows full faithfulness for $p$-complete gauges.
For essential surjectivity, we see that the essential image of $\sigma^*$ contains all $p$-torsion objects, hence all $p^n$-torsion objects by the already proved full faithfulness, and hence all $p$-complete objects.

We therefore wish to prove that the functor
\[
\sigma^* \from\bD(\mathcal{D}^{(m-1,m),\star}_k)\to \bD(\mathcal{D}^{(m,m+1),\star}_k) 
\]
is an equivalence of categories if $m\geq 1$, and induces an equivalence of categories
\[
\bD_{\rm qn}(\mathcal{D}^{(-1,0),\star}_k)\to \bD_{\rm qn}(\mathcal{D}^{(0,1),\star}_k)
\]
if $m=0$. Since both functors are clearly exact, it is enough to prove the corresponding result for static objects.

By Theorems \ref{theo: cdF 1} and \ref{theo: cdF 2} we know that, on each graded piece, the functor $M^i\mapsto \sigma^*M^i$ is an equivalence (either for all objects if $m\geq 1$, or for quasi-nilpotent objects if $m=0$), so it suffices to show that the maps $\partial^{[p^{m}]}\from M^i\to M^{i+1}$ can be uniquely recovered from $\partial^{[p^{m+1}]}\from \sigma^*M^i\to\sigma^*M^{i+1}$. If $m\geq 1$, this can be done in the usual way \cite[Proposition 2.2.4]{Ber00}, since we have the explicit reconstruction
\[
M^i=\left\{ x\in \sigma^*M^i \,\left\vert\, \partial^{\langle j \rangle_{(m)}}x=0\;\; \forall p\nmid j \right. \right\}
\]
of each $M^i$ from $\sigma^*M^{i}$, and the map $\partial^{[p^{m}]}\from M^i\to M^{i+1}$ is then simply the restriction of the map $\partial^{[p^{m+1}]}\from \sigma^*M^i\to\sigma^*M^{i+1}$ to the subspace $M^i\subset \sigma^*M^i$.

If $m=0$, then we use \cite[Theorem 5.8]{GLSQ10} instead to recover each $M^i$ from $\sigma^*M^{i}$. One can then check by an explicit calculation that $\partial^{[p]}\from \sigma^*M^{i}\to \sigma^*M^{i+1}$ restricts to a map $M^i\to M^{i+1}$, where it is given by $(1+g\theta)\partial$ for some $g\in \widehat{\mathcal{D}}^{(-1)}_C$. Since $M^{i+1}$ is quasi-nilpotent, the map $(1+g\theta)$ is invertible on $M^{i+1}$, thus $\partial\from M^i\to M^{i+1}$ can be recovered from $\partial^{[p]}$.
\end{proof}

\section{Symmetric monoidal structure and duality}
\label{sec: symm}

To endow the category of ($F$-)gauges over $k$ with a symmetric monoidal structure, we set $\mathcal{B}^{(m)}\subset \widehat{\mathcal{D}}_C^{(m-1,m),\star}$ to be the subring generated by the images of $C,f,v$, all $\partial^{\langle j\rangle_{(m-1)}}$, and $\partial^{[p^m]}$. The category of condensed graded left $\mathcal{B}^{(m)}$-modules has a natural symmetric monoidal structure $\otimes_{C^\star}$, where each $\partial^{\langle j\rangle_{(m-1)}}$ acts via the standard Leibniz rule
\[ 
\partial^{\langle j\rangle_{(m-1)}}(x\otimes y)=\sum_{i=0}^j \left\{\begin{smallmatrix}j \\ i \end{smallmatrix}\right\}_{(m-1)} \partial^{\langle j-i\rangle_{(m-1)}}x \otimes \partial^{\langle i\rangle_{(m-1)}}y
\]
and $\partial^{[p^m]}$ acts via the graded Leibniz rule
\[
\partial^{[p^m]}(x\otimes y) = \partial^{[p^m]}x \otimes y + f\sum_{j=1}^{p^m-1} \partial^{[p^m-j]}x \otimes \partial^{[j]}y +  x \otimes \partial^{[p^m]}y
\]
Note that this makes sense since $\partial^{[j]}\in \widehat{\mathcal{D}}_C^{(m-1)}$ whenever $1\leq j\leq p^m-1$. Taking the $\infty$-categorical derived functor of $\otimes_{C^\star}$ (essentially appealing to \cite[Proposition 4.1.7.4]{Lur17}) we get a symmetric monoidal structure $\bD(\mathcal{B}^{(m)})^\otimes$ on the derived $\infty$-category of left $\mathcal{B}^{(m)}$-modules, the tensor product for which will be denoted $\otimes^{\mathbf L}_{C^\star}$.

Since $\widehat{\mathcal{D}}_C^{(m-1,m),\star}$ is equal to the $p$-completion of $\mathcal{B}^{(m)}$, the category $\bD(\widehat{\mathcal{D}}_C^{(m-1,m),\star})$ of gauges of level $m$ can be viewed as a reflective subcategory of $\bD(\mathcal{B}^{(m)})$, with derived $p$-completion providing the left adjoint to the inclusion. In particular, $\bD(\widehat{\mathcal{D}}_C^{(m-1,m),\star})$ is a localisation of $\bD(\mathcal{B}^{(m)})$ by \cite[Proposition 5.2.7.12]{Lur09}. We thus define a map $N_1\to N_2$ in $\bD(\mathcal{B}^{(m)})$ to be a $p$-complete equivalence if it becomes an equivalence after derived $p$-completion, or equivalently if the derived $p$-completion of its cofibre is $0$. 

\begin{lemma} For every $M\in \bD(\mathcal{B}^{(m)})$, and every $p$-complete equivalence $N_1\to N_2$, the map $M\otimes_{C^\star}^{\bf L} N_1\to M\otimes_{C^\star}^{\bf L} N_2$ is a $p$-complete equivalence.
\end{lemma}

\begin{proof} Since $M\otimes_{C^\star}^{\mathbf{L}}$ commutes with cofibre sequences, we need to show that if $\widehat{N}=0$, then $\left(M\otimes_{C^\star}^{\bf L}N\right)^\wedge=0$. But vanishing of a derived $p$-complete complex can be checked on its mod $p$ reduction, so we need to show that if $N/^{\bf L}p=0$, then $(M\otimes^{\bf L}_{C^\star} N)/^{\bf L}p=0$, which is clear.
\end{proof} 

Once again applying \cite[Proposition 4.1.7.4]{Lur17}, we obtain a symmetric monoidal structure on $\bD(\widehat{\mathcal{D}}_C^{(m-1,m),\star})$. Concretely, this is given by taking the derived tensor product $\otimes^{\bf L}_{C^\star}$ and then $p$-completing, we will therefore denote this tensor product by $\widehat{\otimes}^{\bf L}_{C^\star}$. Essentially by definition, the gauge $C^\star=C[f,v]/(fv-p)$ is a unit for this tensor product.

\begin{lemma} $\bD(\widehat{\mathcal{D}}_C^{(m-1,m),\star})^\otimes$ is a presentably symmetric monoidal $\infty$-category.
\end{lemma}

\begin{proof} Since $p$-adic completion creates colimits in $\bD(\widehat{\mathcal{D}}_C^{(m-1,m),\star})$, and is (by construction) symmetric monoidal, it is enough to show that $\bD(\mathcal{B}^{(m)})^\otimes$ is presentably symmetric monoidal. Since the underlying category $\bD(\mathcal{B}^{(m)})$ is stable and presentable, it is enough by \cite[Proposition 4.4.2.7]{Lur09} to show that for any given object $M\in \bD(\mathcal{B}^{(m)})$, the functor $M\otimes_{C^\star}^{\bf L}$ commutes with cofibre sequences and direct sums, both of which are clear.
\end{proof}

Of course, we can endow each $\bD(\widehat{\mathcal{D}}_C^{(m)})$ with the structure of a presentably symmetric monoidal $\infty$-category in exactly the same way, and the functors $\widehat{(\,\cdot\,)}^{\pm\infty}$ and $\sigma^*$ are all symmetric monoidal, and all preserve colimits. Combining \cite[Proposition 5.5.3.13]{Lur09} with \cite[Corollary 3.2.2.4]{Lur17} we may therefore define a presentably symmetric monoidal structure $\bD_F(\widehat{\mathcal{D}}_C^{(m-1,m),\star})^\otimes$ on the $\infty$-category $\bD_F(\widehat{\mathcal{D}}_C^{(m-1,m),\star})$ of $F$-gauges by taking the equaliser of the diagram
\[ \bD(\widehat{\mathcal{D}}^{(m-1,m),\star}_C)^\otimes \xrightrightarrows[\sigma^*\widehat{(\,\cdot\,)}^{-\infty}]{\widehat{(\,\cdot\,)}^{+\infty}} \bD(\widehat{\mathcal{D}}^{(m)}_C)^\otimes   \]
In particular, the forgetful functor
\[\bD_F(\widehat{\mathcal{D}}^{(m-1,m),\star}_C) \to \bD(\widehat{\mathcal{D}}^{(m-1,m),\star}_C) \]
is symmetric monoidal (and, of course, commutes with colimits). As a presentably symmetric monoidal $\infty$-category, internal homs exist in $\bD_F(\widehat{\mathcal{D}}^{(m-1,m),\star}_C)^\otimes$, and we will denote them by ${\bf R}{\rm Hom}_{C^\star}(\,\cdot\,,\,\cdot\,)$.

\begin{definition} Let $M$ be an ($F$-)gauge. We let $M^\vee={\bf R}{\rm Hom}_{C^\star}(M,C^\star)$ denote the dual of $M$.
\end{definition}

We can use perfectness of the underlying complex of graded $C^\star$-modules to show that lisse ($F$-)gauges are dualisable, thus the formation of $M^\vee$ for such gauges commutes with the functor which forgets the Frobenius structure. This also shows that the dual of a lisse ($F$-)gauge is again lisse, and that lisse ($F$-)gauges are reflexive in the sense that the double duality map
\[  M\lto (M^\vee)^\vee\]
is an equivalence.

We can give a reasonably explicit description of the dual of a gauge $M$ over $k$, using what we call the Matlis dual of $M$.
We consider the static gauge $E^\star$ given by the diagram 
\[
\ldots  \overset{1}{\underset{p}{\rightleftarrows}} C \overset{1}{\underset{p}{\rightleftarrows}} C  \overset{p}{\underset{1}{\rightleftarrows}} C \overset{p}{\underset{1}{\rightleftarrows}} \ldots
\]
with central term sitting in graded degree $0$. 
The $\widehat{\mathcal{D}}^{(m-1)}_C$-module structure in each degree is the canonical one, and the action of $\partial^{[p^m]}$ is $\frac{1}{c_m}\partial^{\langle p^m\rangle_{(m-1)}}$ in non-negative degrees, and the usual action in negative degrees. Since $\widehat{E}^{\star,\pm\infty}=0$, this gauge has a trivial $F$-structure. 

\begin{definition}
\label{defn: Mat dual}
We define the Matlis dual of a gauge $M$ by
\[
M^\lozenge := {\bf R}{\rm Hom}_{C^\star}(M,E^\star) \in \bD(\widehat{\mathcal{D}}_C^{(m-1,m),\star})
\]
\end{definition}

This has a very natural interpretation as the `degree-wise' dual of $M$. That is, the degree $i$ part of $M^\lozenge$ is the ordinary $C$-linear dual
\[
M^{\lozenge,i}= {\bf R}{\rm Hom}_C(M^{-i},C)=:M^{-i,\vee}
\]
with its natural structure as a (complex of) $\widehat{\mathcal{D}}_C^{(m-1)}$-module(s).
The map $f\from M^{\lozenge,i}\to M^{\lozenge,i+1}$ is given by dualising $f\from M^{-(i+1)}\to M^{-i}$, and the action of $v$ on $M^\lozenge$ is obtained similarly.
The map $\partial ^{[p^m]}\from M^{\lozenge,i}={\bf R}{\rm Hom}_C(M^{-i},C)\to M^{\lozenge,i+1}={\bf R}{\rm Hom}_C(M^{-(i+1)},C)$ is given, roughly speaking, by the `Leibniz rule'
\[
(\partial^{[p^m]}\psi)(x) = 
\partial^{[p^m]}(\psi(fx)) -\psi(\partial^{[p^m]}x)-\sum_{j=1}^{p^m-1}(\partial^{[p^m-j]}\psi)(\partial^{[j]}fx)
\]
for $\psi\in {\bf R}{\rm Hom}_C(M^{-i},C)$ and $x \in M^{-(i+1)}$. This description as the degree-wise dual shows that Matlis duality preserves quasi-lisse gauges, and that quasi-lisse gauges are reflexive for Matlis duality, in the sense that the double duality map
\[
M \to(M^\lozenge)^\lozenge
\]
is an equivalence.
Note that if $M$ is a lisse $F$-gauge, then $\widehat{M}^{\lozenge,\pm\infty}=0$, and so $M^\lozenge$ has a trivial $F$-structure.
Alternatively, we can use the fact that $M$ is dualisable in $\bD_{F}(\widehat{\mathcal{D}}_C^{(m-1,m),\star})$ to define
\[
M^\lozenge = {\bf R}{\rm Hom}_{C^\star}(M,E^\star) \in \bD_{F}(\widehat{\mathcal{D}}_C^{(m-1,m),\star})
\] 
as internal hom in the category of $F$-gauges.
Forgetting the Frobenius structure then recovers the Matlis dual of $M$ as given by Definition \ref{defn: Mat dual}.

To relate the Matlis dual $M^\lozenge$ to the ordinary dual $M^\vee$, we use the natural exact sequence of static $F$-gauges
\begin{equation}
\label{eqn: Matlis}
0 \to C^\star \to  j_+^{(m)}C \oplus j_+^{(m-1)}C \to E^\star \to 0
\end{equation}
The individual terms $j_+^{(m)}C$ and $j_+^{(m-1)}C$ here do not have $F$-structures, but we can put an $F$-structure on their direct sum using the fact that
\[
\left(j_+^{(m)}C \oplus j_+^{(m-1)}C \right)^{\wedge,\pm\infty} = C
\]
For any gauge $M$, we have identifications
\begin{align*}
{\bf R}{\rm Hom}_{C^\star}(M,j_+^{(m)}C) &= j_+^{(m)}{\bf R}{\rm Hom}_{C}(\widehat{M}^{+\infty},C)\\
{\bf R}{\rm Hom}_{C^\star}(M,j_+^{(m-1)}C) &= j_+^{(m-1)}{\bf R}{\rm Hom}_{C}(\widehat{M}^{-\infty},C)
\end{align*}
As before let us write $\widehat{M}^{\pm\infty,\vee}={\bf R}{\rm Hom}_{C}(\widehat{M}^{\pm\infty},C)$ for the $C$-linear dual of $\widehat{M}^{\pm\infty}$, equipped with its natural structure as a $\widehat{\mathcal{D}}_C^{(m)}$ or $\widehat{\mathcal{D}}_C^{(m-1)}$-module.
Then the exact sequence \eqref{eqn: Matlis} gives rise to:

\begin{proposition}
Let $M$ be a gauge of level $m$ over $k$. Then there is a natural fibre sequence
\[
M^\vee \xrightarrow{(f^\infty,v^\infty)} 
j^{(m)}_+(\widehat{M}^{+\infty,\vee})\oplus j^{(m-1)}_+(\widehat{M}^{-\infty,\vee}) \xrightarrow{(f^{\infty,\vee},-v^{\infty,\vee})}
M^\lozenge \overset{+1}{\lto}
\]
\end{proposition}

Since a lisse gauge $M$ satisfies $\widehat{M}^{\lozenge,\pm\infty}=0$, we deduce:

\begin{corollary}
Let $M$ be a lisse gauge of level $m$ over $k$.
Then $M^{\vee,\pm\infty} \simeq M^{\pm\infty,\vee}$.
\end{corollary}

\section{Independence of the lift}\label{sec: indep}

It isn't too hard to use similar ideas to those of \S\ref{sec: cond F desc} to deduce from established independence results for abstract $\mathcal{D}$-modules that the various categories of (condensed) gauges defined here are independent of the choice of lift $(C,\sigma)$. 

First suppose that $m\geq 1$. Then abstract Frobenius descent implies that the level raising module $\sigma^*\widehat{\mathcal{D}}^{(m-1)}_{C'}$ is independent of $\sigma$ (up to canonical isomorphism) as an abstract bimodule, and thus (since it is of finite type) as a condensed bimodule.
This implies that the pullback functor $\sigma^*$ used to form the $\infty$-category of $F$-gauges of level $m$ is independent of the choice of $\sigma$. 

To show independence of the choice of Cohen ring $C$, suppose that we have two such Cohen rings $C_1,C_2$ for $k$. Then there exists an isomorphism $\beta\from C_2\isomto C_1$ inducing the identity on $k$, although since $k$ is not perfect, such an isomorphism is not unique. The basic idea, as in \cite[\S2.1.6]{Ber00}, is to show that the base change bimodule associated to $\beta$ is in fact independent of $\beta$, up to canonical isomorphism.

Indeed, exactly as in \cite[\S5]{Dod22}, we may form the graded $(\widehat{\mathcal{D}}^{(m-1,m),\star}_{C_1},\widehat{\mathcal{D}}^{(m-1,m),\star}_{C_2})$-bimodule 
\[
\widehat{\mathcal{D}}^{(m-1,m),\star}_{\beta} := C_1\otimes_{C_2,\beta} \widehat{\mathcal{D}}^{(m-1,m),\star}_{C_2}
\]
Note that this construction works equally well on the level of condensed bimodules, and recovers the expected abstract bimodule upon taking sections over a point.
This bimodule then induces an equivalence
\[
\bD_{(F)}(\widehat{\mathcal{D}}_{C_2}^{(m-1,m),\star}) \isomto \bD_{(F)}(\widehat{\mathcal{D}}_{C_1}^{(m-1,m),\star})
\]
given by taking the graded tensor product 
\[
M \mapsto \widehat{\mathcal{D}}^{(m-1,m),\star}_{\beta}  \otimes_{ \widehat{\mathcal{D}}^{(m-1,m),\star}_{C_2} } M =  \widehat{\mathcal{D}}^{(m-1,m),\star}_{\beta}  \widehat{\otimes}^{\bf L}_{ \widehat{\mathcal{D}}^{(m-1,m),\star}_{C_2} } M
\]
In fact, this can be upgraded to an equivalence
\begin{equation}
\label{eqn: beta equiv}
\bD_{(F)}(\widehat{\mathcal{D}}_{C_2}^{(m-1,m),\star})^\otimes \to \bD_{(F)}(\widehat{\mathcal{D}}_{C_1}^{(m-1,m),\star})^\otimes
\end{equation}
of symmetric monoidal $\infty$-categories in the obvious way.

We now suppose that we have two such isomorphisms $\beta,\beta'\from C_2\isomto C_1$, both inducing the identity on $k$, and consider each graded piece
\[
\widehat{\mathcal{D}}^{(m-1,m),i}_{\beta} =  C_1\otimes_{C_2,\beta} \widehat{\mathcal{D}}^{(m-1,m),i}_{C_2}
\]
as a $(\widehat{\mathcal{D}}^{(m-1)}_{C_1},\widehat{\mathcal{D}}^{(m-1)}_{C_2})$-bimodule.
If $m\geq 1$ and $p>2$, or if $m\geq 2$, then  case a) of \cite[Proposition 2.1.5]{Ber00} provides a canonical isomorphism
\[
\tau^i_{\beta,\beta'}\from \widehat{\mathcal{D}}^{(m-1,m),i}_{\beta'}\isomto \widehat{\mathcal{D}}^{(m-1,m),i}_{\beta}
\]
of abstract $(\widehat{\mathcal{D}}^{(m-1)}_{C_1},\widehat{\mathcal{D}}^{(m-1)}_{C_2})$-bimodules, which is necessarily continuous and therefore upgrades to an isomorphism of condensed bimodules. 

\begin{proposition}
The isomorphisms $\tau^i_{\beta,\beta'}$ are compatible with the actions of $f,v$ and $\partial^{[p^m]}$.
They therefore induce an isomorphism
\[
\tau_{\beta,\beta'}\from \widehat{\mathcal{D}}^{(m-1,m),\star}_{\beta'}\isomto \widehat{\mathcal{D}}^{(m-1,m),\star}_{\beta}
\]
of condensed $(\widehat{\mathcal{D}}^{(m-1,m),\star}_{C_1},\widehat{\mathcal{D}}^{(m-1,m),\star}_{C_2})$-bimodules. 
\end{proposition}

\begin{proof}
As usual, it suffices to prove compatibility on the underlying abstract modules. In this case, using the definition of $\widehat{\mathcal{D}}^{(m-1,m),i}_{C_2}$ as a submodule of $\widehat{\mathcal{D}}^{(m)}_{C_2}$, the proposition follows from applying \cite[Proposition 2.1.5]{Ber00} instead to $\widehat{\mathcal{D}}^{(m)}_{C_2}$ as a module over itself, and then noting that Berthelot's isomorphism $\tau_{\beta,\beta'}$ is compatible with restriction of scalars along ${\mathcal{D}}^{(m-1)}\to \mathcal{D}^{(m)}$. 
\end{proof}

Again, it follows from case a) of \cite[Proposition 2.1.5]{Ber00} that we have $\tau_{\beta,\beta}=\mathrm{id}$ and $\tau_{\beta,\beta''}=\tau_{\beta,\beta'}\circ \tau_{\beta',\beta''}$, and so the object $\widehat{\mathcal{D}}^{(m-1,m),\star}_{\beta}$ is independent of any choices, up to canonical isomorphism. It follows that the equivalence \eqref{eqn: beta equiv} is independent of the choice of $\beta$.

\begin{corollary}
Suppose that either $m\geq 2$, or $m\geq 1$ and $p>2$. Then the $\infty$-category $\bD_{(F)}(\widehat{\mathcal{D}}^{(m-1,m),\star}_C)$ of ($F$-)gauges of level $m$ is independent of the choice of $(C,\sigma)$ up to canonical equivalence. 
\end{corollary}

It is tedious, but not difficult, to extend this to an equivalence of symmetric monoidal $\infty$-categories.  
If $p=2$ and $m=1$ then the argument breaks down since the $m$-PD nilpotence hypothesis in case a) of \cite[Proposition 2.1.5]{Ber00} need not hold.
But we do still get a similar result for quasi-nilpotent ($F$-)gauges, since we can simply use case b) of \cite[Proposition 2.1.5]{Ber00} instead. 
We can also extend these independence results to quasi-nilpotent gauges of level $0$, using Frobenius descent for gauges. 

\begin{corollary}
For any prime $p$, and any level $m\geq 0$, the $\infty$-category $\bD_{(F),\rm qn}(\widehat{\mathcal{D}}^{(m-1,m),\star}_C)$ of quasi-nilpotent ($F$-)gauges of level $m$ is independent of $(C,\sigma)$, up to canonical equivalence.
\end{corollary}

Again, this can be upgraded to include symmetric monoidal structures.

\section{Nygaard cohomology, \texorpdfstring{de\thinspace Rham}{de Rham} cohomology and syntomic cohomology}
\label{sec: cohom}

A gauge or $F$-gauge of level $0$ over $k$ has various natural cohomology groups associated to it, all of which are given by taking ${\bf R}{\rm Hom}(C^\star,\,\cdot\,)$ in a suitable category mapped to from $\bD_{(F)}(\widehat{\mathcal{D}}^{(-1,0),\star}_C)$.
In fact, by taking internal homs in the category of condensed sets, these cohomology groups will all come with a natural condensed structure.
In this section we will define these cohomology groups, and give explicit complexes calculating them.
Of course, gauges of higher levels will also have the same associated cohomology groups, but these will be harder to compute, and, by Frobenius descent, will give nothing new.

\subsection{Condensed cohomology of $\bm{F}$-gauges}

We let $\bD(\Z_p)$ (resp.\ $\bD(W)$) denote the full subcategory of the derived $\infty$-category of condensed $\Z_p$-modules (resp.\ $W$-modules) consisting of derived $p$-complete complexes. We similarly let $\bD(\F_p)$ (resp.\ $\bD(\F)$) denote the derived $\infty$-category of condensed $\F_p$-modules (resp.\ $\F$-modules). As before, we will generally understand $W$, $\Z_p$, $\F$ and $\F_p$ to refer to the condensed version of these rings, and so will refer to, for example, a complex of $W$-modules to mean a complex of condensed $W$-modules.

\begin{definition}
Let $M\in  \bD(\widehat{\mathcal{D}}^{(-1,0),\star}_C)$. We define the Nygaard cohomology of $M$ to be
\[ {\bf R}\Gamma_{\EU N}(k,M):={\bf R}{\rm Hom}_{\widehat{\mathcal{D}}^{(-1,0),\star}_C} (C^\star,M)  \in \bD(W)  \]
This is the (condensed) $W$-linear complex of graded homomorphisms $C^\star\to M$. \end{definition}

\begin{definition}
We similarly define the level $0$ and level $-1$ de\thinspace Rham cohomology of $M$ to be the condensed $W$-linear complexes
\begin{align*}
{\bf R}\Gamma_{\dR,(0)}(k,M) &:={\bf R}{\rm Hom}_{\widehat{\mathcal{D}}^{(0)}_C}(C,\widehat{M}^{+\infty})  \\
{\bf R}\Gamma_{\dR,(-1)}(k,M) &:={\bf R}{\rm Hom}_{\widehat{\mathcal{D}}^{(-1)}_C}(C,\widehat{M}^{-\infty})  
\end{align*}
These are given by homomorphisms $C^\star\to M$ after applying the functors $\widehat{(\,\cdot\,)}^{+\infty}$ and $\widehat{(\,\cdot\,)}^{-\infty}$.
\end{definition}

We will also use the notation ${\bf R}\Gamma_{\dR,(m)}(k,\,\cdot\,)$ for $p$-complete complexes of $\widehat{\mathcal{D}}^{(m)}_C$-modules, thus by definition we have
\begin{align*}
{\bf R}\Gamma_{\dR,(0)}(k,M) &={\bf R}\Gamma_{\dR,(0)}(k,\widehat{M}^{+\infty})  \\
{\bf R}\Gamma_{\dR,(-1)}(k,M) &={\bf R}\Gamma_{\dR,(-1)}(k,\widehat{M}^{-\infty})
\end{align*}
whenever $M$ is a gauge. Note that adjunction gives
\begin{equation}
\label{eqn: Nyg dR}
{\bf R}\Gamma_{\dR,(m)}(k,M)\isomto{\bf R}\Gamma_{\EU N}(k,j^{(m)}_+M)
\end{equation}
and thus de\thinspace Rham cohomology can be viewed as a special case of Nygaard cohomology. Syntomic cohomology is then just absolute cohomology in the category of $F$-gauges.

\begin{definition}
Let $M\in  \bD_F(\widehat{\mathcal{D}}^{(-1,0),\star}_C)$. We define the syntomic cohomology of $M$ to be
\[
{\bf R}\Gamma_{\rm Syn}(k,M):={\bf R}{\rm Hom}_{ \bD_F(\widehat{\mathcal{D}}^{(-1,0),\star}_C)}(C^\star,M) \in \bD(\Z_p)
\]
This is the condensed $\Z_p$-linear complex of homomorphisms $C^\star\to M$ in  $\bD_F(\widehat{\mathcal{D}}^{(-1,0),\star}_C)$. 
\end{definition}

Since $C^\star$ is the unit object for the symmetric monoidal structure on $F$-gauges, all of these cohomology theories are lax monoidal.

Note that $f$ and $v$ induce maps
\begin{align*} f^\infty\from  {\bf R}\Gamma_{\EU N}(k,M) &\to{\bf R}\Gamma_{\dR,(0)}(k,M)\\ 
 v^\infty\from  {\bf R}\Gamma_{\EU N}(k,M) &\to{\bf R}\Gamma_{\dR,(-1)}(k,M)  \end{align*}
By the definition of $\bD_F(\widehat{\mathcal{D}}^{(-1,0),\star}_C)$ as an equaliser, there is then a natural fibre sequence
\[
{\bf R}\Gamma_{\rm Syn}(k,M) \lto
{\bf R}\Gamma_{\EU N}(k,M)\xrightarrow{f^\infty-\phi^\flat v^\infty} 
 {\bf R}\Gamma_{\dR,(0)}(k,M) \overset{+1}{\lto}
 \]
in $\bD(\Z_p)$, for any $F$-gauge $M$. Here $\phi^\flat$ is the composite
\[
{\bf R}{\rm Hom}_{\widehat{\mathcal{D}}^{(-1)}_C}(C,\widehat{M}^{-\infty}) \lto {\bf R}{\rm Hom}_{\widehat{\mathcal{D}}^{(0)}_C}(C,\sigma^*\widehat{M}^{-\infty}) \isomto {\bf R}{\rm Hom}_{\widehat{\mathcal{D}}^{(0)}_C}(C,\widehat{M}^{+\infty})
\]
where the first map is the semilinear map coming from Frobenius base change, and the second map is the $W$-linear equivalence induced by the $F$-structure $\phi\from \sigma^*\widehat{M}^{-\infty}\isomto \widehat{M}^{+\infty}$. If $M$ is quasi-nilpotent, then the first map (and hence the composite $\phi^\flat$) is an equivalence.

\subsection{Computing cohomology} To compute these various cohomology groups explicitly, we need to find a resolution of $C^\star$ by free graded $\widehat{\mathcal{D}}^{(-1,0),\star}_C$-modules. To do so, we choose a uniformiser lift $t\in C(*)$, and consider the corresponding differential operator $\partial$. The kernel of the surjective map 
\begin{align*}
\widehat{\mathcal{D}}^{(-1,0),\star}_C &\to C^\star= C[f,v]/(fv-p) \\
P&\mapsto P\cdot 1
\end{align*}
is then the left ideal generated by the degree $1$ element $\partial$. We therefore obtain the standard (co-ordinate free) presentation
\[
0 \lto \widehat{\mathcal{D}}^{(-1,0),\star}_C(1) \otimes_C \widehat{\ca{T}}_{C} \lto \widehat{\mathcal{D}}^{(-1,0),\star}_C \lto C^\star \lto 0
\]
of $C^\star$ as a graded left $\widehat{\mathcal{D}}^{(-1,0),\star}_C$-module.
Here $\widehat{\ca{T}}_{C}={\rm Hom}_C(\widehat{\Omega}^1_C,C)$ is the $p$-complete tangent module of $C$, and the first map in this sequence sends $P\otimes \partial$ to $P\partial$. 

\begin{proposition}
\label{prop: cohom explicit 1}
Let $M$ be a gauge of level $0$ over $k$.
Then the Nygaard and de\thinspace Rham cohomologies of $M$ are given by the totalisations of the following double complexes
\begin{align*}
{\bf R}\Gamma_{\EU N}(k,M)&= \left[ M^0 \overset{\nabla}{\lto} \widehat{\Omega}^1_C\otimes_C M^1\right] \\
{\bf R}\Gamma_{\dR,(0)}(k,M)&= \left[ \widehat{M}^{+\infty} \overset{\nabla}{\lto}  \widehat{\Omega}^1_C\otimes_C \widehat{M}^{+\infty} \right] \\
{\bf R}\Gamma_{\dR,(-1)}(k,M)&= \left[ \widehat{M}^{-\infty} \overset{\nabla^{(-1)}}{\lto}  \widehat{\Omega}^1_C\otimes_C \widehat{M}^{-\infty}\right]
\end{align*}
Here $\nabla^{(-1)}$ is the map defined by $\nabla^{(-1)}(x)={\d}t\otimes \theta(x)$.
\end{proposition}

For a general gauge $M$, the maps 
\begin{align*}
f^\infty &\from  {\bf R}\Gamma_{\EU N}(k,M)\to{\bf R}\Gamma_{\dR,(0)}(k,M) \\
v^\infty &\from  {\bf R}\Gamma_{\EU N}(k,M)\to {\bf R}\Gamma_{\dR,(-1)}(k,M)
\end{align*}
are the obvious ones, and if $M$ has an $F$-structure, then the map
\[
\phi^\flat \from  {\bf R}\Gamma_{\dR,(-1)}(k,M) \to{\bf R}\Gamma_{\dR,(0)}(k,M) 
\]
is given by the diagram
\[ 
\xymatrix{
\widehat{M}^{-\infty} \ar[r]^-{\nabla^{(-1)}}\ar[d]_{\phi^\flat} & \widehat{\Omega}^1_C\otimes_C\widehat{M}^{-\infty} \ar[d]^{\frac{\sigma}{p}\otimes \phi^\flat} \\
\widehat{M}^{+\infty} \ar[r]^-{\nabla}& \widehat{\Omega}^1_C\otimes_C\widehat{M}^{+\infty}
}
\]

\begin{corollary}
\label{cor: syn coh explicit}
Let $M$ be an $F$-gauge of level $0$ over $k$. Then ${\bf R}\Gamma_{\rm Syn}(k,M)$ is represented by the totalisation of the triple complex 
\[ 
\xymatrix{
 M^0 \ar[r]^-{\nabla}  \ar[d]_-{f^\infty-\phi^\flat v^\infty} &
\widehat{\Omega}^1_C\otimes_C M^1   \ar[d]^-{1\otimes f^\infty- \frac{\sigma}{p}\otimes \phi^{\flat}v^\infty} \\
 \widehat{M}^{+\infty} \ar[r]^-{\nabla} &
 \widehat{\Omega}^1_C\otimes_C \widehat{M}^{+\infty}
 }
 \]
\end{corollary}

\section{Virtual \texorpdfstring{$F$}{F}-crystals, Dieudonn\'e modules, and \texorpdfstring{$F$}{F}-gauges}
\label{sec: vfc etc}

In this section we compare the category of $F$-gauges over $k$ to other categories of integral $p$-adic coefficients, in particular virtual $F$-crystals, $p$-divisible groups and finite flat group schemes (of $p$-power order). In each case we construct a fully faithful functor into the category of $F$-gauges, and show that syntomic cohomology as defined above computes the expected cohomology groups. 

\subsection{Virtual $\bm{F}$-crystals and $\bm{F}$-gauges}

As in \cite{FJ21}, one important source of $F$-gauges over $k$ is virtual $F$-crystals.

\begin{definition}
A virtual $F$-crystal over $k$ is a finite free $C$-module $\mathcal{M}$, together with a connection and a horizontal isomorphism
\[
\varphi\from \sigma^*\mathcal{M}[1/p]\isomto \mathcal{M}[1/p]
\]
The category of such objects will be denoted ${\bf vF}\text{-}{\bf Cris}_k$. 
\end{definition} 

Perhaps the most natural construction of an $F$-gauge from a virtual $F$-crystal $\mathcal{M}$ occurs at level $1$. Given a virtual $F$-crystal $\mathcal{M}$, we define a gauge $M:=\fr{g}^{(1)}\mathcal{M}$ of level $1$ using the standard Nygaard filtration
\[
M^{-i}
:={\rm F}^{i}_\EU{N}\sigma^*\mathcal{M}
:= \left\{ x\in \sigma^*\mathcal{M} \,\left\vert\, \varphi(x)\in p^{i}\mathcal{M}\right.\right\} 
\]
The maps $f\from M^i\to M^{i+1}$ are the natural inclusions, and the maps $v\from M^i\to M^{i-1}$ are multiplication by $p$.
An easy check shows that $\partial(M^i)\subset M^i$ and $\frac{\partial^p}{p!}M^i\subset M^{i+1}$, and that these maps give rise to the structure of a gauge of level $1$ on $M=\bigoplus_i M^i$. 
There is a canonical identification $M^{+\infty}=\sigma^*\mathcal{M}$ as $\widehat{\mathcal{D}}^{(1)}_C$-modules. Moreover, the various maps $p^i\varphi\from M^i\to \mathcal{M}$ fit together to give an isomorphism $M^{-\infty}\isomto \mathcal{M}$ of $\widehat{\mathcal{D}}^{(0)}_C$-modules.
The gauge $M$ therefore inherits an $F$-structure from the identity map $\sigma^*\mathcal{M}\isomto \sigma^*\mathcal{M}$, and we get a functor
\begin{equation}
\label{eqn: g1}
\fr{g}^{(1)}\from {\bf vF}\text{-}{\bf Cris}_k \to \bD_{F,{\rm lisse}}(\widehat{\mathcal{D}}^{(0,1),\star}_C)
\end{equation}
from virtual $F$-crystals to $F$-gauges over $k$.
This functor is readily checked to be fully faithful, however, as in \cite{FJ21}, the essential image is smaller than one might initially guess. 

\begin{definition}
Let $M$ be a static gauge (of any level $m$) over $k$.
We say that $M$ is rigid if the maps
\[
f\from M^i/v \to M^{i+1}/v,\qquad v\from M^i/f\to M^{i-1}/f
\]
are injective, for all $i$.
\end{definition}

This is a slight abuse of the same terminology from \cite{FJ21}. The proof of the following result is, however, essentially identical to that of \cite[Theorem 3.2.9]{FJ21}.

\begin{theorem} Let $M$ be an $F$-gauge of level $1$ over $k$. Then $M$ is in the essential image of \eqref{eqn: g1} if and only if it is static, lisse, $p$-torsion free and rigid.
\end{theorem}

To calculate the cohomology of such gauges, we note that by Theorem \ref{theo: F desc gauge}, there exists a unique functor
\[
\fr{g}^{(0)}\from {\bf vF}\text{-}{\bf Cris}_k \to \bD_{F,{\rm lisse}}(\widehat{\mathcal{D}}^{(-1,0),\star}_C)
\]
taking values in $F$-gauges of level $0$, such that $\sigma^*\fr{g}^{(0)}\isomto \fr{g}^{(1)}$. To describe the functor $\fr{g}^{(0)}$, we note that $\ca{M}$, being quasi-nilpotent by Lemma \ref{lemma: qn nilp}, itself has a canonical Frobenius descent. Thus we get a finite free $C$-module $\ca{N}$, equipped with a `$p$-connection' (in other words, a compatible $\widehat{\mathcal{D}}^{(-1)}_C$-module structure) and a horizontal isomorphism
\[
\varphi_{-1} \from \sigma^*\ca{N}[1/p]\isomto \ca{N}[1/p] 
\]
such that $\ca{M}=\sigma^*\ca{N}$ and $\varphi=\sigma^*\varphi_{-1}$. Then 
\[
(\fr{g}^{(0)}\ca{M})^{-i} = \left\{ x\in \mathcal{M} \,\left\vert\, \varphi_{-1}(x)\in p^{i}\mathcal{N}\right.\right\}
\]
and the actions of $f,v,\theta$ and $\partial$ can be described exactly as in the case of $\fr{g}^{(1)}$. To relate this to the more familiar Nygaard filtration on $\ca{M}$, we let $\varphi^\flat$ denote the composite
\[ 
\ca{M} \xrightarrow{\iota(x)=1\otimes x} \sigma^*\ca{M} \overset{\varphi}{\lto} \ca{M}[1/p]
\]
and define 
\[
{\rm F}^{i}_\EU{N}\mathcal{M}:=\left\{ x\in \mathcal{M} \,\left\vert\, \varphi^\flat(x)\in p^{i}\mathcal{M}\right.\right\}
\]
The fact that $\varphi=\sigma^*\varphi_{-1}$ means that the diagram
\[ 
\xymatrix{ 
\ca{M}[1/p] \ar[r]^\iota \ar[d]_{\varphi_{-1}} & \sigma^*\ca{M}[1/p] \ar[d]^{\varphi} \\
\ca{N}[1/p] \ar[r]^-\iota & \ca{M}[1/p]=\sigma^*\ca{N}[1/p]
}
\]
commutes, and finite flatness of $\sigma$ then gives an equality $(\fr{g}^{(0)}\ca{M})^{-i}={\rm F}^{i}_\EU{N}\mathcal{M}$ for all $i$. This means that we can compute syntomic cohomology of $\fr{g}^{(0)}\ca{M}$ in the expected manner, using the divided Frobenius maps
\[
\varphi^\flat_i := \frac{\varphi^\flat}{p^i} \from {\rm F}^i_{\EU N}\mathcal{M}\to \mathcal{M}
\]
Indeed, Corollary \ref{cor: syn coh explicit} now shows:

\begin{theorem}
\label{theo: syn cohom complex 2}
Let $\mathcal{M}$ be a virtual $F$-crystal over $k$.
Then ${\bf R}\Gamma_{\rm Syn}(k,\fr{g}^{(0)}\mathcal{M})$ is represented by the totalisation of the double complex 
\[
\xymatrix{
{\rm F}^0_\EU{N}\mathcal{M}\ar[r]^-\nabla \ar[d]_{1-\varphi^\flat_0} & \widehat{\Omega}^1_C\otimes_C{\rm F}^{-1}_\EU{N}\mathcal{M} \ar[d]^{1-\frac{\sigma}{p}\otimes \varphi^\flat_{-1}} \\
\mathcal{M} \ar[r]^-{\nabla} & \widehat{\Omega}^1_C\otimes_C\mathcal{M}
}
\]
\end{theorem}

\subsection{Dieudonn\'e modules and \texorpdfstring{$\bm{F}$}{F}-gauges}

As explained in the introduction, local Tate duality for syntomic cohomology should be seen as an extension of local Tate duality for finite flat group schemes proved in \cite{Sha64}. The purpose of this section is to show how to recover the cohomology of finite flat group schemes over $k$ as a special case of the syntomic cohomology of $F$-gauges, using the above results on the cohomology of virtual crystals.

Recall that a Dieudonn\'e module over $C$ is a finitely generated $C$-module $M$, together with a connection $\nabla\from M\to \widehat{\Omega}^1_C\otimes_C M$ and horizontal maps $F\from \sigma^*M\to M$ and $V\from M\to \sigma^*M$ such that $FV=VF=p$, the category of such objects is denoted ${\bf DM}_C$. Of course, we may view all such objects with their canonical condensed structure, and thus view ${\bf DM}_C$ as a category of condensed $C$-modules together with $F,V$ and $\nabla$ as above. This category is then the target of the Dieudonn\'e module functor, defined for both $p$-divisible groups and finite flat group schemes over $k$. Indeed, if we let ${\bf BT}_k$ denote the category of $p$-divisible groups over $k$, and ${\bf FG}_k$ the category of (commutative) finite flat group schemes over $k$ of $p$-power order, then evaluating the contravariant Dieudonn\'e crystal functor of \cite[\S3]{BBM82} on the pro-PD-thickening $C\twoheadrightarrow k$ gives rise to functors
\begin{align*}   \mathbb{M}^*&\from {\bf BT}^{\rm op}_k \to {\bf DM}_C \\
 \mathbb{M}^*&\from {\bf FG}_k^{\rm op} \to {\bf DM}_C
\end{align*}
with the former landing inside the full subcategory of Dieudonn\'e modules which are finite free over $C$, and the latter inside the full subcategory of those which are of finite length over $C$. In both cases, we define $\mathbb{M}(G):=\mathbb{M}^*(G^D)$ where $(\,\cdot\,)^D$ denotes the Cartier dual functor. These functors are then compatible in the expected way, in that  
\[ \mathbb{M}(G)=\lim{n}\mathbb{M}(G[p^n]) \]
whenever $G$ is a $p$-divisible group, and any exact sequence
\[ 0\to H\to G_1\to G_2 \to 0\]
where $H$ is finite flat and $G_1,G_2$ are $p$-divisible gives rise to an exact sequence
\[ 0\to \mathbb{M}(G_1)\to \mathbb{M}(G_2)\to \mathbb{M}(H) \to 0 \]
Both Dieudonn\'e module functors are fully faithful \cite[Th\'eor\`eme 4.1.1]{BM90}, and that for $p$-divisible groups is an equivalence onto the full subcategory of objects which are finite free over $C$ \cite[Main Theorem 1]{dJ95b} (see also the proof of \cite[Proposition 2.4.8]{dJ95b} for why quasi-nilpotence of the connection is automatic). The Dieudonn\'e module functor for finite flat group schemes is \emph{not} an equivalence, even if one imposes quasi-nilpotence of the connection, as we shall see below.

As with virtual $F$-crystals, the natural functor from Dieudonn\'e modules to gauges lives at level $1$. For concreteness, we fix a uniformiser lift $t\in C(*)$, and use the corresponding differential operators $\partial$ and $\partial^{[p]}$.

\begin{definition} \label{defn: DM gaugelike}
We say that a Dieudonn\'e module $M\in {\bf DM}_C$ is gauge-like if $(p-1)!F\partial^{[p]}V=\partial^p$ as maps $M\to M$. We denote the full subcategory of such objects by ${\bf DM}_C^{\fr g}$.
\end{definition}

Note that this condition is automatic if $M$ is torsion-free, since we always have
\[ p!F\partial^{[p]}V= F\partial^pV=FV\partial^p = p\partial^p \]
We also note that it is stable under both subobjects and quotients, and doesn't depend on the choice of uniformiser lift $t$. There do, however, exist torsion Dieudonn\'e modules which are not gauge-like. In fact, there are even quasi-nilpotent examples.

\begin{example}
\label{exa: non gauge-like}
Let $M=ke_1\oplus ke_2$ be free of rank $2$ over $k$, with connection given by $\nabla(e_1)=\d t\otimes t^{p-1}e_2$, $\nabla(e_2)=0$.
Set $F=V=0$.
Then $M$ is a Dieudonn\'e module over $C$.
Since $F=V=0$ we see immediately that $(p-1)!F\partial^{[p]}V=0$, however, $\partial^p(e_1)=\partial^{p-1}(t^{p-1})e_2=-e_2\neq 0$. 
Hence $M$ is not gauge-like, although it is quasi-nilpotent.
\end{example}

Nevertheless, all objects coming from finite flat group schemes will be gauge-like.
 
\begin{lemma}
\label{lemma: psirig}
Let $G$ be either a finite flat group scheme or a $p$-divisible group over $k$. Then $\mathbb{M}(G)$ is gauge-like.
\end{lemma}

\begin{proof}
If $G$ is a $p$-divisible group, then $\mathbb{M}(G)$ is torsion-free, and hence gauge-like.
If $H$ is a finite flat group scheme, then by \cite[Th\'eor\`eme 3.1.1]{BBM82} there exists an injective map $H\hookrightarrow G$ with $G$ a $p$-divisible group.
This induces a surjective map $\mathbb{M}(G/H)\to \mathbb{M}(H)$, where $G/H$ is $p$-divisible, and the fact that $\mathbb{M}(H)$ is gauge-like follows from the fact that $\mathbb{M}(G/H)$ is.
\end{proof}

The utility of Definition \ref{defn: DM gaugelike} is that it is precisely the condition that allows a Dieudonn\'e module to be converted into a (static) $F$-gauge of level $1$ over $k$.
That is, given a Dieudonn\'e module $M$ over $C$, we consider the graded $\widehat{\mathcal{D}}^{(0)}_C[f,v]/(fv-p)$-module represented by the diagram
\[
\ldots  \overset{p}{\underset{1}{\rightleftarrows}} M \overset{p}{\underset{1}{\rightleftarrows}} M \overset{V}{\underset{F}{\rightleftarrows}} \sigma^*M \overset{1}{\underset{p}{\rightleftarrows}} \sigma^*M \overset{1}{\underset{p}{\rightleftarrows}} \ldots 
\]
where $M$ is in graded degree $0$ and $\sigma^*M$ in graded degree $1$ (note that $V$ plays the role of $f$ and $F$ that of $v$).
We can define maps $\partial^{[p]}\from \sigma^*M\to \sigma^*M$ using the canonical $\widehat{\mathcal{D}}^{(1)}_C$-module structure on the Frobenius pullback of the $\widehat{\mathcal{D}}^{(0)}_C$-module $M$, and maps $\partial^{[p]}\from M\to M$ as $\frac{1}{(p-1)!}\partial^p$. We can define an $F$-structure on this object by taking the identity map $\sigma^*M\to\sigma^*M$. 

The final piece of the $F$-gauge structure is then a map
\[ \partial^{[p]} \from M\to \sigma^*M, \]
which we define to be $\partial^{[p]}V$, where $\partial^{[p]}\from \sigma^*M\to \sigma^*M$ is given by the above natural $\widehat{\mathcal{D}}^{(1)}_C$-module structure on $\sigma^*M$. That the relations necessary to have a $\widehat{\mathcal{D}}^{(0,1),\star}_C$-module structure are satisfied then amounts to the identity
\[ (p-1)!F\partial^{[p]}V = \partial^p \from M\to M \]
which is exactly the statement that $M$ is gauge-like. This construction therefore gives rise to a fully faithful functor
\[ {\bf DM}_C^{\fr g}\lto \bD_{F,{\rm lisse}}(\widehat{\mathcal{D}}^{(0,1),\star}_C) \]
landing in the subcategory of static, lisse $F$-gauges of level 1 over $k$. In fact, one easily checks that this is an equivalence onto the full subcategory of static lisse $F$-gauges which are concentrated in the interval $[0,1]$, in the sense that $f\from M^i\to M^{i+1}$ is an isomorphism whenever $i\geq 1$, and $v\from M^i\to M^{i-1}$ is an isomorphism whenever $i\leq 0$. By Lemma \ref{lemma: psirig}, we may apply this construction to the Dieudonn\'e modules arising from $p$-divisible groups or finite flat group schemes over $k$.

\begin{definition}
We define $\fr{g}^{(1)}\from {\bf FG}_k\to \bD_{F,{\rm lisse}}(\widehat{\mathcal{D}}^{(0,1),\star}_C)$ and $\fr{g}^{(1)}\from {\bf BT}_k\to \bD_{F,{\rm lisse}}(\widehat{\mathcal{D}}^{(0,1),\star}_C)$ as the composite functors
\begin{align*}
{\bf FG}_k& \overset{\mathbb{M}}{\lto}{\bf DM}_C^{\fr g}\lto \bD_{F,{\rm lisse}}(\widehat{\mathcal{D}}^{(0,1),\star}_C)  \\
 {\bf BT}_k &\overset{\mathbb{M}}{\lto}{\bf DM}_C^{\fr g}\lto \bD_{F,{\rm lisse}}(\widehat{\mathcal{D}}^{(0,1),\star}_C)
 \end{align*}
\end{definition}

In both cases, $\fr{g}^{(1)}$ is a fully faithful, covariant functor.
As in the case of Dieudonn\'e modules, we have $\fr{g}^{(1)}G=\lim{n}\fr{g}^{(1)}G[p^n]$ whenever $G$ is a $p$-divisible group, and any exact sequence 
\[
0 \to H \to G_1\to G_2 \to 0
\]
where $H$ is finite and $G_1,G_2$ are $p$-divisible gives rise to an exact sequence
\[
0 \to \fr{g}^{(1)}G_1\to \fr{g}^{(1)}G_2\to \fr{g}^{(1)}H \to 0
\]
of static $F$-gauges.
For the standard finite flat group schemes over $k$, $\fr{g}^{(1)}$ behaves in the expected manner.

\begin{example}
\label{exa: standard ffgs}
We consider the finite flat group schemes $\Z/p^n$, $\mu_{p^n}$ and $\alpha_p$ over $k$. Then $\fr{g}^{(1)}\Z/p^n=C^\star/p^n$, $\fr{g}^{(1)}\mu_{p^n}=C^\star(1)/p^n$, and $\fr{g}^{(1)}\alpha_p$ is the $F$-gauge given by
\[
\ldots  \overset{0}{\underset{1}{\rightleftarrows}} k \overset{0}{\underset{1}{\rightleftarrows}} k  \overset{0}{\underset{0}{\rightleftarrows}} k \overset{1}{\underset{0}{\rightleftarrows}} \ldots 
\]
where the central term sits in graded degree $0$.
\end{example}

Again, by Theorem \ref{theo: F desc gauge}, there exists a functor $\fr{g}^{(0)}$ from either $p$-divisible groups or finite flat group schemes to static, lisse, $F$-gauges of level $0$, uniquely determined by the requirement that $\fr{g}^{(1)}=\sigma^*\fr{g}^{(0)}$. 
The computations of $\fr{g}^{(0)}$ in the standard cases of $\Z/p^n$, $\mu_{p^n}$ and $\alpha_p$ are exactly as in Example \ref{exa: standard ffgs}, in particular we have $\fr{g}^{(0)}\Z/p^n=C^\star/p^n$ and $\fr{g}^{(0)}\mu_{p^n}=C^\star(1)/p^n$.

The flat cohomology of a $p$-divisible group $G$ over $k$ can be calculated in terms of its covariant Dieudonn\'e module $M:=\mathbb{M}(G)$. As for a virtual $F$-crystal, we let $\iota\from M\to \sigma^*M$ be the map $x\mapsto 1\otimes x$, and set $F^\flat:=F\iota\from M\to M$ to be the induced semilinear map. The Hodge filtration on $M$ can then be defined by the rule
\[
{\rm Fil}^iM:=\left\{ m\in M \,\left\vert\, F^\flat(m) \in p^iM \right.\right\} \]
As with virtual $F$-crystals, there exists a $\sigma$-semilinear map
\[ \frac{F^\flat}{p}\from {\rm Fil}^1M \to M
\]
defined by $\frac{F^\flat}{p}=\frac{1}{p}F\iota$.
Then \cite[Theorem 1.1]{TV24} states that, if $p>2$, then the flat cohomology complex ${\bf R}\lim{n}{\bf R}\Gamma(k,G[p^n]) $ can be computed by (the abstract complex of $\Z_p$-modules underlying) the double complex
\begin{equation*}
\xymatrix{
{\rm Fil}^1M \ar[r]^-\nabla \ar[d]_{1-\frac{F^\flat}{p}}& \widehat{\Omega}^1_C\otimes_C M \ar[d]^{1-\frac{\sigma}{p}\otimes F^\flat}
\\ M \ar[r]^-\nabla &\widehat{\Omega}^1_C\otimes_C M 
}
\end{equation*}
In fact, it is explained in \cite[p.662, final paragraph]{BKK26} that the restriction $p>2$ is unnecessary. Applying Theorem \ref{theo: syn cohom complex 2} to the virtual $F$-crystal $(M,p^{-1}F)$, we therefore have:

\begin{theorem}
Let $G$ be a $p$-divisible group over $k$. Then there is a natural equivalence
\[
{\bf R}\lim{n}{\bf R}\Gamma(k,G[p^n]) \isomto {\bf R}\Gamma_{\rm Syn}(k,\fr{g}^{(0)}G)(*)
\]
of complexes of abstract $\Z_p$-modules.
\end{theorem}

The left hand side here should be considered as a kind of continuous cohomology ${\bf R}\Gamma_{\rm cont}(k,T_pG)$. Via the usual trick of embedding a finite group scheme in a $p$-divisible group, we therefore obtain:

\begin{corollary}
\label{cor: syn coh flat}Let $H$ be a finite flat group scheme over $k$ of $p$-power order. Then there is a natural equivalence
\[
{\bf R}\Gamma(k,H) \isomto {\bf R}\Gamma_{\rm Syn}(k,\fr{g}^{(0)}H)(*)
\]
of complexes of abstract $\Z_p$-modules.
\end{corollary}

At least after taking cohomology groups, we will upgrade this to a topological isomorphism in \S\ref{sec: shatz comp} below.

\section{Pontryagin duality and the trace map}
\label{sec: PD TM}

Ultimately, local Tate duality in syntomic cohomology will boil down to Pontryagin duality for finitely generated $k$-vector spaces. In this section, we interpret this duality in the language of condensed mathematics, and carefully construct the trace map as a morphism
\[ {\bf R}\Gamma_{\rm Syn}(k,C^\star(1)[2]) \to \Z_p \]
of complexes of condensed $\Z_p$-modules. From now on, we will only be interested in ($F$-)gauges of level $0$, which we will therefore refer to simply as ($F$-)gauges.

\subsection{Pontryagin duality for $\bm{k}$-modules}\label{sec: pd k}

The functor $N\mapsto N(*)$ is an equivalence of categories between finitely generated condensed $k$-modules and finite dimensional abstract $k(*)$-vector spaces. We may therefore view any finitely generated $k(*)$-module canonically as a condensed $k$-module. In particular, we obtain the invertible $k$-module $\Omega^1_k$ of absolute differentials on $k$. 
 
 Any choice of uniformiser $t\in k(*)$ induces an isomorphism $k \cong \F\lser{t} \isomto \Omega^1_k$ given by $\alpha \mapsto \alpha\,\dlog t$. We may therefore define the usual residue map ${\Res}\from \Omega^1_k\to \F$ by sending $\sum_i a_i t^i\,\dlog t$ to $a_{0}$. The results of \cite[Chapter II]{Ser88} show that this is independent of the choice of uniformiser $t$. 
 
 \begin{definition} We define the trace map ${\Tr}\from \Omega^1_k\to \F_p$ as the composite of ${\Res}$ with the usual Galois-theoretic trace map ${\Tr}_{\F/\F_p}\from \F\to \F_p$. 
 \end{definition}
 
 Pontryagin duality for finitely generated $k$-modules then takes the following form. 
 
 \begin{proposition} \label{prop: pd k} Let $N$ be a perfect complex of $k$-modules, with $k$-linear dual complex $N^\vee:={\bf R}{\rm Hom}_k(N,k)$. Then the map \[ \left(\Omega^1_k\otimes_k N^\vee \right)\otimes^{\bf L}_k N  \overset{\rm ev}{\lto} \Omega^1_k \overset{\Tr}{\lto} \F_p \]
 induces an equivalence
 \[ \Omega^1_k  \otimes_k N^\vee  \isomto {\bf R}{\rm Hom}_{\F_p}(N,\F_p) \]
in ${\bD}(\F_p)$.
 \end{proposition}
 
\begin{proof}
We immediately reduce to the case $N=k$ and $\F=\F_p$. Choosing a uniformiser $t$, the given pairing can be identified with the map
\begin{align*}
\F_p\lser{t} \times \F_p\lser{t} &\to \F_p \\
\left(\sum_i a_it^i ,\sum_i b_it^i\right)&\mapsto \sum_i a_ib_{-i}
\end{align*}
We therefore reduce to showing that the pairing
\begin{align*}
\bigoplus_\N \F_p  \times \prod_\N \F_p &\to \F_p \\
\left(a_i, b_i \right)&\mapsto \sum_i a_ib_i
\end{align*}
induces equivalences
\begin{align*}
\prod_{\N} \F_p &\isomto {\bf R}{\rm Hom}_{\F_p}(\bigoplus_\N \F_p,\F_p)\\
\bigoplus_{\N} \F_p &\isomto {\bf R}{\rm Hom}_{\F_p}(\prod_\N \F_p,\F_p)
\end{align*}
The first of these follows more or less immediately from the definitions, using the fact that countable products of (light) condensed abelian groups are exact. 

For the second, it follows from \cite[proof of Theorem 3.2.4(9)]{RC26} that the right hand side is concentrated in degree $0$, and a straightforward calculation shows that the given map is an isomorphism in this degree.  \end{proof}

\subsection{The trace map in syntomic cohomology}

We now extend the trace map from \S\ref{sec: pd k} to a map
\[
\Tr\from {\bf R}\Gamma_{\rm Syn}(k,C^\star(1)[2] )\to \Z_p
\]
in $\bD(\Z_p)$.
To do so, we note that Corollary \ref{cor: syn coh explicit} shows that ${\bf R}\Gamma_{\rm Syn}(k,C^\star(1))$ is given explicitly as the totalisation of the double complex
\[
\xymatrix{
pC \ar[r]^-{{\d}}\ar[d]_{1-\frac{\sigma}{p}} & \widehat{\Omega}^1_C \ar[d]^{1-\frac{\sigma}{p}} \\ 
C \ar[r]^-{\d} &  \widehat{\Omega}^1_C
}
\]
We define the residue map
\[
\Res \from \widehat{\Omega}^{1}_C\to W
\]
in the usual way, taking $\sum_i \alpha_it^i \dlog t$ to $\alpha_0$ whenever $t\in C(*)$ is a uniformiser lift.
To prove independence of $\Res$ from the choice of $t$, we will use the following simple lemma.

\begin{lemma} 
\label{lemma: check dense}
Let $\psi,\xi\from C\to W$ be homomorphisms of condensed $W$-modules, and let $t\in C(*)$ be a uniformiser lift.
If $\psi(t^i)=\xi(t^i)$ for all $i\in \Z$, then $\psi=\xi$. 
\end{lemma}

\begin{proof}
Since $W$ is $p$-adically complete, it is enough to show that $\psi=\xi$ as maps $C/p^n\to W/p^n$, for all $n$.
Next, writing $C/p^n = \colim{m} t^{-m} W/p^n\pow{t}$ it is enough to show that $\psi=\xi$ as maps $t^{-m} W/p^n\pow{t}\to W/p^n$, for all $m$ and $n$.
But now $t^{-m} W/p^n\pow{t}$ and $W/p^n$ are both light profinite sets, and so it suffices to observe that $W/p^n$ is Hausdorff and that $t^{-m} W/p^n[t] \subset t^{-m} W/p^n\pow{t}$ is dense.
\end{proof}

Note that since $\sigma$ is an automorphism of $W$, we can (and will) also apply the lemma to $\sigma$-linear homomorphisms $C\to W$.

\begin{lemma}
The map $\Res$ is independent of $t$.
\end{lemma} 

\begin{proof}
By Lemma \ref{lemma: check dense} it is enough to show that if $s\in C(*)$ is another such choice of uniformiser lift, then
\[
\Res(s^i\dlog s)=\begin{cases} 1& i=0 \\ 0& i\neq 0\end{cases}
\]
If $i\neq 0$ then this is easy: we simply calculate  
\[
is^i\,\dlog s=\d s^i = t\partial(s^i)\,\dlog t
\]
and then observe that $t\partial(s^i)$ clearly has no constant term, and that $W$ is torsion-free.

If $i=0$, we write $s=tu_0u_1$ where $u_0\in W\pow{t}^\times$ and $u_1\in 1+pC$.
Then 
\[
\log u_1 = -\sum_{n\geq 1} \frac{(1-u_1)^n}{n}
\]
converges in $C$, and so we find that
\[
\dlog s = \dlog t + tu_0^{-1}\partial u_0\cdot \dlog t + \d( \log u_1)
\]
and the second and third terms visibly have zero residue in $t$.
\end{proof}

\begin{lemma}
The map $\Tr_{W/\Z_p}\circ \Res\from \widehat{\Omega}^{1}_C\to \Z_p $ vanishes on the images of both $\d\from C\to\widehat{\Omega}^{1}_C$ and $1-\frac{\sigma}{p}\from  \widehat{\Omega}^1_C\to \widehat{\Omega}^{1}_C$.
\end{lemma}

\begin{proof}
That $\Tr_{W/\Z_p}\circ \Res\circ \d=0$ is immediate, since we already have $\Res\circ \d=0$ by inspection.
We next show that ${\rm Res}$ intertwines $\frac{\sigma}{p}$ with the usual Witt vector Frobenius $\sigma$ on $W$. Indeed, by Lemma \ref{lemma: check dense}, it is enough to show that ${\rm Res}(\frac{\sigma}{p}(t^i\,\dlog t))=0$ for $i\neq 0$, and that ${\rm Res}(\frac{\sigma}{p}(\dlog t))=1$. For the former, we calculate  
\[
i\sigma(t^i\,\dlog t) = \sigma (\d t^i) = \d\sigma(t^i)
\]
which shows that $ip\cdot {\rm Res}(\frac{\sigma}{p}(t^i\,\dlog t))=0$. Thus ${\rm Res}(\frac{\sigma}{p}(t^i\,\dlog t))=0$ since $W$ is torsion-free. For the latter, we write $\sigma(t)=t^pu$ for some $u\in 1+pC$. As before,
\[ \frac{\log u}{p} = -\sum_{n\geq 1} \frac{p^{n-1}\left(\frac{1-u}{p}\right)^n}{n} \]
converges in $C$, and so
\[ \frac{\sigma}{p}(\dlog t)=  \dlog t +\d(\frac{\log u}{p})\]
which clearly has residue $1$. To complete the proof, we now observe that, for any $\omega\in \widehat{\Omega}^1_C$,
\[ {\rm Res}((1-\tfrac{\sigma}{p})\omega) = (1-\sigma){\rm Res}(\omega) \]
and ${\rm Tr}_{W/\Z_p}$ vanishes on the image of $1-\sigma$ by Galois theory for $W/\Z_p$.
\end{proof}

\begin{definition} We define the trace map
\[ \Tr\from {\bf R}\Gamma_{\rm Syn}(k,C^\star(1)[2]) \to \Z_p \]
to be the composite
\[ {\bf R}\Gamma_{\rm Syn}(k,C^\star(1)[2])\lto \rH^2_{\rm Syn}(k,C^\star(1)) \xrightarrow{ \Tr_{W/\Z_p}\circ \Res} \Z_p \]
\end{definition}

\subsection{Pontryagin duality for $\bm{C}$-modules, and Matlis duality revisited}

We can now combine Pontryagin duality for $k$-modules with the residue map constructed over $C$ to deduce a form of Pontryagin duality for perfect complexes of $C$-modules (or, equivalently, perfect complexes of $C(*)$-modules). For any such complex $M$, we write $M^\vee:={\bf R}{\rm Hom}_C(M,C)$ for its (condensed) $C$-linear dual, which is again a complex of the same form. 

\begin{proposition}
\label{prop: pd C}
Let $M$ be a perfect complex of $C$-modules. Then the map
\[ \left( \widehat{\Omega}^1_C \otimes_C M^\vee\right) \otimes_C^{\bf L} M\overset{\rm ev}{\lto} \widehat{\Omega}^1_C \xrightarrow{\Tr_{W/\Z_p}\circ \Res} \Z_p  \]
induces an equivalence
\[ \widehat{\Omega}^1_C \otimes_C M^\vee \isomto {\bf R}{\rm Hom}_{\Z_p}(M,\Z_p) \]
in ${\bD}(\Z_p)$.
\end{proposition}

\begin{proof} This reduces via a $p$-completeness argument to Proposition \ref{prop: pd k}.
\end{proof}

We denote by $(\,\cdot\,)^*={\bf R}{\rm Hom}_{\Z_p}(\,\cdot\,,\Z_p)$ the duality functor for (condensed) $\Z_p$-modules. One nice consequence of Proposition \ref{prop: pd C} is then the reflexivity of Nygaard, de\thinspace Rham and syntomic cohomology of gauges, subject to appropriate finiteness conditions. 

\begin{corollary}
Let $M\in \bD(\widehat{\mathcal{D}}_C^{(-1,0),\star})$ be a quasi-lisse gauge over $k$. 
Then the double duality map
\[
{\bf R}\Gamma_{\EU N}(k,M)\to {\bf R}\Gamma_{\EU N}(k,M)^{**}
\]
is an equivalence.
If $M^{\pm\infty}$ is a perfect complex of $C$-modules, then the same is true for de\thinspace Rham cohomology, and if $M$ is a lisse $F$-gauge, then the same is true for syntomic cohomology.
\end{corollary}

\begin{proof}
Since de\thinspace Rham cohomology is a special case of Nygaard cohomology, by \eqref{eqn: Nyg dR}, and syntomic cohomology sits in a fibre sequence with Nygaard and de\thinspace Rham cohomology, it suffices to treat the case of Nygaard cohomology.
We can then use Proposition \ref{prop: cohom explicit 1} to compute the Nygaard cohomology of $M$ as the double complex
\[
M^0 \overset{\nabla}{\lto} \widehat{\Omega}^1_C\otimes_C M^1
\]
We now simply note that $M^0$ and $M^1$ are perfect complexes over $C$, and hence reflexive by Proposition \ref{prop: pd C}. 
\end{proof}
 
We can also use Pontryagin duality to give an alternative interpretation of the Matlis dual functor $M^\lozenge$ introduced in \S\ref{sec: symm}.
Indeed, if $M$ is a gauge, then we can form its degree-wise $\Z_p$-linear dual $M^*$ via the rule
\[
M^{*,i}={\bf R}{\rm Hom}_{\Z_p}(M^{-i},\Z_p)
\]
The graded left $\widehat{\mathcal{D}}^{(-1,0),\star}_C$-module structure on the complex $M$ gives rise to a graded right $\widehat{\mathcal{D}}^{(-1,0),\star}_C$-module structure on $M^*$, simply by taking the $\Z_p$-dual of the action of each element $P\in \widehat{\mathcal{D}}^{(-1,0),\star}_C$.
On the other hand, the side switching module $\widehat{\Omega}^1_C$ makes $\widehat{\Omega}^1_C \otimes_C M^\lozenge$ into a complex of graded right $\widehat{\mathcal{D}}^{(-1,0),\star}_C$-modules as in \cite[\S6]{Dod22}. If $M$ is quasi-lisse, then taking the map from Proposition \ref{prop: pd C} degree-wise gives a natural map
\begin{equation}
\label{eqn: Mat linear}
\widehat{\Omega}^1_C \otimes_C M^\lozenge \to M^*
\end{equation}
of complexes of graded $C$-modules.

\begin{proposition}
\label{prop: star-diamond}
Let $M$ be a quasi-lisse gauge over $k$.
Then the map \eqref{eqn: Mat linear} is $\widehat{\mathcal{D}}^{(-1,0),\star}_C$-linear.
In particular, there is a natural equivalence
\[
M^\lozenge \isomto M^* \otimes_C \widehat{\Omega}^{1,\vee}_C
\]
of quasi-lisse gauges.
\end{proposition}

\begin{proof}
We need to show that for any $P\in \widehat{\mathcal{D}}^{(-1,0),\star}_C$, the map
\[
\xymatrix{
\widehat{\Omega}^1_C \otimes_C M^\lozenge \ar[d]_{\cdot P}\ar[r] & M^* \ar[d]^{(P\cdot )^*} \\
\widehat{\Omega}^1_C \otimes_C M^\lozenge \ar[r] & M^* 
}
\]
commutes, for which it is sufficient to show that 
\[
{\rm Res}(\omega(Px)) = {\rm Res}((\omega P)x)
\]
for any $\omega \in \widehat{\Omega}^1_C  \otimes_C M^\lozenge$, $x\in M$. We already know $C$-linearity, so it suffices to treat the three cases $P=f,v,\partial$.
The first two are straightforward. 

For the last, if we have $\omega\in \widehat{\Omega}^1_C \otimes_C M^{\lozenge,i}=\widehat{\Omega}^1_C \otimes_C M^{-i,\vee}$, and $x\in M^{-(i+1)}$, then the Leibniz rule for right $\mathcal{D}$-modules gives
\[
(\omega(fx))\partial =  (\omega\partial)x- \omega(\partial x)  
\]
and so 
\[
{\rm Res}((\omega\partial)x) - {\rm Res}( \omega(\partial x)) = {\rm Res}((\omega(fx))\partial) =0
\]
as required. 
\end{proof}

\section{Local Tate duality}\label{sec: core}

We now have all the tools we need to state and prove our main theorem. If $N$ is any complex of (condensed) $\Z_p$-modules, recall the notation
\[
N^*:={\bf R}{\rm Hom}_{\Z_p}(N,\Z_p)
\]
for the $\Z_p$-linear dual of $N$. 

\begin{theorem}[Local Tate duality]
\label{theo: main1} Let $M$ be a lisse $F$-gauge over $k$.
Then the evaluation pairing
\[
M\widehat{\otimes}^{\bf L}_{C^\star} M^\vee(1) \to C^\star(1)
\]
composed with the trace map 
\[
\Tr\from {\bf R}\Gamma_{\rm Syn}(k,C^\star(1)[2]) \to \Z_p
\]
induces an equivalence
\[
{\bf R}\Gamma_{\rm Syn}(k,M^\vee(1)[2])  \isomto {\bf R}\Gamma_{\rm Syn}(k,M)^*
\]
in $\bD(\Z_p)$.
\end{theorem}

\subsection{Nygaard duality}

Before proving Theorem \ref{theo: main1}, we will need an auxiliary duality result on the Nygaard cohomology of quasi-lisse gauges. Recall from \S\ref{sec: symm} above the gauge $E^\star$, represented by the diagram
\[
\ldots  \overset{1}{\underset{p}{\rightleftarrows}} C \overset{1}{\underset{p}{\rightleftarrows}} C  \overset{p}{\underset{1}{\rightleftarrows}} C \overset{p}{\underset{1}{\rightleftarrows}} \ldots
\]
We have $\widehat{E}^{\star,\pm\infty}=0$, and so ${\bf R}\Gamma_{\dR,(m)}(k,E^\star)=0$ for $m=-1,0$. Hence the morphism $E^\star\to C^\star[1]$ of $F$-gauges coming from the exact sequence \eqref{eqn: Matlis}  induces a map
\[ {\bf R}\Gamma_{\EU N}(k,E^\star(1)[1])= {\bf R}\Gamma_{\rm Syn}(k,E^\star(1)[1]) \lto {\bf R}\Gamma_{\rm Syn}(k,C^\star(1)[2])\overset{\rm Tr}{\lto} \Z_p \]
which we also call the trace map.

\begin{proposition} \label{prop: N dual} Let $M$ be a quasi-lisse gauge over $k$. Then the pairing
\[
{\bf R}\Gamma_{\EU N}(k,M)\widehat{\otimes}^{\bf L}_{\Z_p}{\bf R}\Gamma_{\EU N}(k,M^\lozenge(1)[1]) \overset{{\rm ev}}{\lto} {\bf R}\Gamma_{\EU N}(k,E^\star(1)[1]) \overset{\rm Tr}{\lto} \Z_p 
\]
induces an equivalence
\[
{\bf R}\Gamma_{\EU N}(k,M^\lozenge(1)[1]) \isomto {\bf R}\Gamma_{\EU N}(k,M)^*
\]
in $\bD(\Z_p)$.
\end{proposition}

\begin{proof}
We may assume that $M$ is static, in which case each graded piece is a finitely generated $C$-module. 
Note that by Propositions \ref{prop: cohom explicit 1} and \ref{prop: star-diamond}, the $\Z_p$-linear dual of 
\[
{\bf R}\Gamma_{\EU N}(k,M) = \left[ M^0 \overset{\nabla}{\lto} \widehat{\Omega}^1_C\otimes_C M^1 \right]
\]
is computed by the totalisation of the complex
\[
\left[ M^{\lozenge,-1} \overset{\nabla}{\lto} \widehat{\Omega}^1_C\otimes_C M^{\lozenge,0} \right]
\]
with $M^{\lozenge,-1}$ sitting in horizontal cohomological degree $-1$, and the implicit pairing being given by the trace map.
Again, by Proposition \ref{prop: cohom explicit 1}, we therefore see that the trace map induces an equivalence
\[ 
{\bf R}\Gamma_{\EU N}(k,M)^*\isomto  {\bf R}\Gamma_{\EU N}(k,M^\lozenge(1)[1])
\]
as claimed.
\end{proof}

\subsection{Proof of local Tate duality} \label{sec: final proof}

We use the explicit description of $M^\vee$ from \S\ref{sec: symm} via the fibre sequence
\[
M^\vee  \lto
j^{(0)}_+(M^{+\infty,\vee}) \oplus j^{(-1)}_+(M^{-\infty,\vee}) \lto
M^{\lozenge} \overset{+1}{\lto}
\]
of gauges.
Twisting, taking Nygaard cohomology, and using \eqref{eqn: Nyg dR}, gives a fibre sequence
\begin{equation*}
{\bf R}\Gamma_{\EU N}(k,M^\vee(1)) \xrightarrow{(f^\infty,v^\infty)} 
{\bf R}\Gamma_{\dR,(0)}(k,M^\vee) \oplus {\bf R}\Gamma_{\dR,(-1)}(k,M^\vee)\xrightarrow{(f^{\infty,\vee},-v^{\infty,\vee})}
{\bf R}\Gamma_{\EU N}(k,M^\lozenge (1)) \overset{+1}{\lto}
\end{equation*}
We now consider the `diagonal' map
\[
{\bf R}\Gamma_{\dR,(0)}(k,M^\vee) \xrightarrow{({\rm id},\phi^{\vee,\flat,-1})} {\bf R}\Gamma_{\dR,(0)}(k,M^\vee)\oplus {\bf R}\Gamma_{\dR,(-1)}(k,M^\vee) 
\]
where $\phi^\vee$ is the Frobenius structure on $M^\vee$.
Pulling back the previous fibre sequence along this map then gives a fibre sequence
\[
{\bf R}\Gamma_{\rm Syn}(k,M^\vee(1)) \lto
{\bf R}\Gamma_{\dR,(0)}(k,M^\vee)\xrightarrow{f^{\infty,\vee}-v^{\infty,\vee}\phi^{\vee,\flat,-1}}
{\bf R}\Gamma_{\EU N}(k,M^\lozenge (1)) \overset{+1}{\lto}
\]
Combining \eqref{eqn: Nyg dR} with Proposition \ref{prop: N dual} we see that the trace map induces an equivalence between this and the fibre sequence
\[
{\bf R}\Gamma_{\rm Syn}(k,M^\vee(1)) \lto
{\bf R}\Gamma_{\dR,(0)}(k,M[1])^* \xrightarrow{(f^\infty-\phi^\flat v^\infty)^*}
{\bf R}\Gamma_{\EU N}(k,M[1])^* \overset{+1}{\lto}
\]
using the fact $f^{\infty,*}=f^{\infty,\vee}$, $v^{\infty,*}=v^{\infty,\vee}$ and $\phi^{\flat,*}=\phi^{\vee,\flat,-1}$, where $\phi$ is the Frobenius on $M$.
Rotating this last fibre sequence then shows that
\[
{\bf R}\Gamma_{\rm Syn}(k,M^\vee(1)) \isomto {\bf R}\Gamma_{\rm Syn}(k,M[2])^*
\]
is an equivalence as claimed.

\subsection{Duality for $\bm{p}$-divisible groups and finite flat group schemes}
\label{sec: shatz comp}

We now show that the duality given by Theorem \ref{theo: main1} restricts to Shatz's duality \cite{Sha64} for finite flat group schemes over $k$ of $p$-power order.
First of all, one easily checks that if $G$ is a $p$-divisible group, then there is a natural isomorphism
\[
\fr{g}^{(0)}G^D \cong (\fr{g}^{(0)}G)^\vee(1)
\]
of $F$-gauges.
As usual, the Tate twist on the right hand side here comes from the fact that Cartier duality is given by homomorphisms into $\G_m$. This in turn gives an isomorphism
\[
\fr{g}^{(0)}H^D \cong (\fr{g}^{(0)}H)^\vee(1)[1]
\]
whenever $H$ is a finite flat group scheme, where the shift on the right hand side now comes from the fact that $\fr{g}^{(0)}H$ is killed by a power of $p$.
Thus, to show that the duality pairing from Theorem \ref{theo: main1} recovers Shatz's duality pairing, it is enough to show that the diagram of abstract groups
\begin{equation*}
\xymatrix{
\rH^2_{\rm Syn}(k,C^\star(1)/p^n)(*) \ar[d]_{\rm Tr} &  \ar[l]_-{\cong} \rH^2(k,\mu_{p^n}) \ar[r] &  \rH^2(k,\G_m) \ar[d]^{\rm inv} \\ 
\Z/p^n \ar[rr]^{p^{-n}} & & \Q/\Z
}
\end{equation*}
commutes, where the top left hand horizontal map is that given by Corollary \ref{cor: syn coh flat}.
The Frobenius lift $\sigma$ gives a map $C\to W(k)$, which induces an isomorphism
\[
\rH^2_{\rm Syn}(k,C^\star(1)/p^n)(*) = \frac{\Omega^1_{C/p^n}}{{\d}(C/p^n)+(1-\frac{\sigma}{p})(\Omega^1_{C/p^n})}\isomto \frac{W_n\Omega^1_k}{{\d}V^{n-1}k+(1-F)W_n\Omega^1_k}
\]
identifying ${\rm Tr}$ with Kato's residue map $\kappa^2_{p^n}$ as defined in \cite[\S2.4]{KM26}.
That this coincides with the classical invariant map was explained in \cite[Example 2.8]{KM26}. 

This then enables us to show that the isomorphism
\begin{equation}
\label{eqn: flat cohom groups syn}
\rH^*(k,H) \isomto \rH^*_{\rm Syn}(k,\fr{g}^{(0)}H)(*)
\end{equation}
from Corollary \ref{cor: syn coh flat} is in fact a homeomorphism, where the right-hand side is equipped with the topology induced by the condensed structure $\rH^*_{\rm Syn}(k,\fr{g}^{(0)}H)$ (that is, the finest topology such that any map $S\to \rH^*_{\rm Syn}(k,\fr{g}^{(0)}H)(*)$ coming from a section of $\rH^*_{\rm Syn}(k,\fr{g}^{(0)}H)$ on a light profinite set $S$ is continuous).
Indeed, both sides are locally profinite abelian groups, and $p^N$-torsion for large enough $N$.
By combining Theorem \ref{theo: main1} with \cite[Theorem on p.411]{Sha64}, we see that \eqref{eqn: flat cohom groups syn} identifies the Pontryagin duals ${\rm Hom}_{\rm cts}(\,\cdot\,,\Q_p/\Z_p)$ of both sides as subgroups of the abstract dual ${\rm Hom}(\,\cdot\,,\Q_p/\Z_p)$.
Thus \eqref{eqn: flat cohom groups syn} must be a homeomorphism by \cite{Hew63}, and we deduce that Theorem \ref{theo: main1} does indeed generalise Shatz's duality result for finite flat group schemes over $k$ of $p$-power order.

\section{Lagrangian refinement of local Tate duality}\label{sec: lagrange}

Another key feature of classical local Tate duality is the fact that the subgroup of unramified classes is self-orthogonal under the duality pairing. For $p$-adic coefficients in mixed characteristic, this is precisely \cite[Corollary 6.4.6]{Bha22}, and for flat cohomology in equicharacteristic the appropriate version of this result is \cite[Corollary III.7.2]{Mil86}. In this section we will prove the analogue of this fact for syntomic cohomology over function fields, as we have defined it using Dodd's category of $F$-gauges.

\subsection{$\bm{F}$-gauges over $\bm{k^\circ}$}

Let $k^\circ$ denote the ring of integers of $k$, and let $C^\circ$ be a $p$-adically complete subring of $C$, stable under $\sigma$, with $C^\circ/p=k^\circ$ (the existence of such a $C^\circ$ is an implicit additional assumption on $(C,\sigma)$, but we can always choose $C$ and $\sigma$ so that such a $C^\circ$ exists). We view both $k^\circ$ and $C^\circ$ as condensed rings in the natural way. Concretely, a choice of uniformiser lift $t\in C^\circ(*)$ induces an isomorphism $C^\circ\cong W\pow{t}$ where $W$, as above, is the Witt ring of the residue field of $k$. Then we have the graded subring
\[
\widehat{\mathcal{D}}^{(-1,0),\star}_{C^\circ} \subset  \widehat{\mathcal{D}}^{(-1,0),\star}_{C}
\]
which is the graded $p$-completion of the subring generated by $C^\circ$, $f$, $v$ and $\partial$.
We may use this to define ($F$-)gauges over $k^\circ$ exactly as over $k$, as well as their Nygaard, de\thinspace Rham and syntomic cohomology groups, using the unit gauge $C^{\circ,\star}:=C^\circ[f,v]/(fv-p)$ over $k^\circ$.
Moreover these cohomologies are all computed by exactly the same formulae as over $k$: the analogues of Proposition \ref{prop: cohom explicit 1}, \eqref{eqn: Nyg dR} and Corollary \ref{cor: syn coh explicit} all hold, with exactly the same proofs. 

\begin{definition}
A(n $F$-)gauge over $k^\circ$ is said to be lisse if it is perfect as a complex of graded $C^{\circ,\star}$-modules.
It is said to be quasi-lisse if it is bounded, and each graded piece is a perfect complex over $C^\circ$. 
\end{definition}

These are rather strong finiteness conditions, for example, they will never be satisfied for non-trivial $F$-gauges obtained by applying the forgetful functor along $ \widehat{\mathcal{D}}^{(-1,0),\star}_{C^\circ} \to  \widehat{\mathcal{D}}^{(-1,0),\star}_{C}$.
It does, however, seem reasonable to expect that $F$-gauges arising as the cohomology of smooth and proper $k^\circ$-schemes should be lisse.
We can at least see that $F$-gauges arising from virtual $F$-crystals, $p$-divisible groups, or finite flat group schemes over $k^\circ$ via functors analogous to those defined in \S\ref{sec: vfc etc} are lisse, and that the analogue of Theorem \ref{theo: syn cohom complex 2} holds, with the same proof.
The base change of a lisse $F$-gauge to $k$ will be lisse, and we may form the dual $M^{\circ,\vee}:={\bf R}{\rm Hom}_{C^{\circ,\star}}(M^\circ,C^{\circ,\star})$ of a lisse $F$-gauge over $k^\circ$ exactly as over $k$.
This dual will again be lisse, and its formation commutes with base change to $k$.

\begin{remark}
\label{rem: tors gauge}
We note that a lisse $F$-gauge over $k^\circ$ can have both $p$-torsion and $t$-torsion. For example, we can consider the gauge
\[
\ldots  \overset{0}{\underset{1}{\rightleftarrows}} k^\circ/(t) \overset{0}{\underset{1}{\rightleftarrows}} k^\circ/(t)  \overset{0}{\underset{0}{\rightleftarrows}} k^\circ/(t^{p}) \overset{1}{\underset{0}{\rightleftarrows}} \ldots
\]
with central term sitting in graded degree $0$, where each term has its canonical $\mathcal{D}^{(-1)}_{k^\circ}$-module structure, and where $\partial$ acts trivially in non-positive degrees, and in the canonical way in positive degrees.
The Frobenius structure is the canonical isomorphism $\sigma^*(k^\circ/(t)) \isomto k^\circ/(t^{p})$. 
\end{remark}

We see that ${\bf R}\Gamma_{\rm Syn}(k^\circ,C^{\circ,\star}(1))$ is computed by the expected complex
\[
\xymatrix{
pC^{\circ} \ar[r]^-{\d}\ar[d]_{1-\frac{\sigma}{p}} & \widehat{\Omega}^1_{C^{\circ}} \ar[d]^{1-\frac{\sigma}{p}} \\ 
C^{\circ} \ar[r]^-{\d} &  \widehat{\Omega}^{1}_{C^{\circ}}
}
\]
and therefore the composite map
\[
\rH^2_{\rm Syn}(k^\circ,C^{\circ,\star}(1)) \to \rH^2_{\rm Syn}(k,C^{\star}(1)) \overset{\rm Tr}{\lto} \Z_p
\]
vanishes. Hence, for any lisse $F$-gauge $M^\circ$ over $k^\circ$, the induced pairing
\[
{\bf R}\Gamma_{\rm Syn}(k^\circ,M^\circ) \widehat{\otimes}_{\Z_p}^{\bf L} {\bf R}\Gamma_{\rm Syn}(k^\circ,M^{\circ,\vee}(1)[2]) \to \Z_p
\]
is just the zero map.
The duality pairing for the base change $M$ of $M^\circ$ to $k$ then gives rise to a commutative square
\begin{equation}
\label{eqn: fibre square}
\xymatrix{
{\bf R}\Gamma_{\rm Syn}(k^\circ,M^\circ) \ar[r] \ar[d] &  {\bf R}\Gamma_{\rm Syn}(k,M) \ar[d] \\ 
0 \ar[r] & {\bf R}\Gamma_{\rm Syn}(k^\circ,M^{\circ,\vee}(1)[2])^*
}
\end{equation}

\begin{theorem}[Lagrangian refinement of local Tate duality]
\label{theo: lagrange}
The square \eqref{eqn: fibre square} is a fibre square.
That is, for any lisse $F$-gauge
 $M^\circ$ over $k^\circ$, with base change $M$ to $k$, the sequence
\begin{equation}
\label{eqn: fibre}
{\bf R}\Gamma_{\rm Syn}(k^\circ,M^\circ) \lto
 {\bf R}\Gamma_{\rm Syn}(k,M)\lto
  {\bf R}\Gamma_{\rm Syn}(k^\circ,M^{\circ,\vee}(1)[2])^* \overset{+1}{\lto}
  \end{equation}
 is a fibre sequence.
\end{theorem}

\subsection{Proof of Theorem \ref{theo: lagrange}}

For a gauge $M^\circ$ over $k^\circ$, we can form the Matlis dual
\[
M^{\circ,\lozenge}:=\bigoplus_i {\bf R}{\rm Hom}_{C^\circ}(M^{\circ,-i},C^\circ)
\]
exactly as we did over $k$ in \S\ref{sec: symm}, and we have the analogous fibre sequence
\[
M^{\circ,\vee} \lto
 j^{(0)}_+(\widehat{M}^{\circ,+\infty,\vee})\oplus j^{(-1)}_+(\widehat{M}^{\circ,-\infty,\vee}) \lto
M^{\circ,\lozenge}\overset{+1}{\lto}
\]
Just as our proof of local Tate duality ultimately reduced to Pontryagin duality for finitely generated $k$-modules, the proof of Theorem \ref{theo: lagrange} will ultimately boil down to Pontryagin duality for $k^\circ$-modules, via the following lemma.

\begin{lemma}
\label{lemma: pont dual orthog}
Let $M^\circ$ be a perfect complex of $C^\circ$-modules, with $C^\circ$-linear dual complex $M^{\circ,\vee}:={\bf R}{\rm Hom}_{C^\circ}(M^\circ,C^\circ)$, and base change $M:=M^\circ\otimes^{\bf L}_{C^\circ}C$ to $C$.
Then the equivalence
\[
 M\isomto \left(\widehat{\Omega}^1_{C} \otimes_C M^\vee\right)^*
\]
from Proposition \ref{prop: pd C} induces a fibre sequence
\[
M^\circ \lto M \lto \left(\widehat{\Omega}^1_{C^\circ} \otimes_{C^\circ}  M^{\circ,\vee}\right)^*\overset{+1}{\lto}
\]
\end{lemma}

\begin{proof}
We may replace $M^\circ$ first by its reduction modulo $p$, and then by its cohomology groups. 
It is then straightforward to reduce to the case $\F=\F_p$, and we can treat separately the case where $M^\circ$ is a finite free $k^\circ\cong \F_p\pow{t}$-module, and where $M^\circ$ is $t^\infty$-torsion. We can reduce the former to the case $M^\circ=\F_p\pow{t}$, and the latter to the case $M^\circ=k^\circ/(t)=\F_p$. 

If $M^\circ=k^\circ$, then the essential claim is that, under the duality pairing
\[
\F_p\lser{t} \times \F_p\lser{t}\, \dlog t \to \F_p,\qquad (\alpha,\beta\,\dlog t)\mapsto {\rm Res}(\alpha\beta \,\dlog t)
\]
the orthogonal complement of $\F_p\pow{t}\subset \F_p\lser{t}$ is $\F_p\pow{t}\, {\d t}\subset \F_p\lser{t}\, \dlog t$, which is a routine calculation. If $M^\circ=\F_p$, then the claim is that 
\[
{\bf R}{\rm Hom}_{\F_p\pow{t}}(\F_p,\Omega^1_{\F_p\pow{t}}) \isomto \F_p[-1]
\]
which is, again, a routine calculation.
\end{proof}

If $M^\circ$ is quasi-lisse, we see that its base change $M$ to $k$ satisfies $(M^\lozenge)^\lozenge =M$, and hence Nygaard duality for $M$ (Proposition \ref{prop: N dual}) gives an equivalence
\[
{\bf R}\Gamma_{\EU N}(k,M) \isomto {\bf R}\Gamma_{\EU N}(k,M^\lozenge(1)[1])^*
\]

\begin{proposition}
\label{prop: N lagrange}
Let $M^\circ$ be a quasi-lisse gauge over $k^\circ$, with base change $M$ to $k$.
Then the equivalence
\[
{\bf R}\Gamma_{\EU N}(k,M) \isomto {\bf R}\Gamma_{\EU N}(k,M^\lozenge(1)[1])^*
\]
of Proposition \ref{prop: N dual} induces a fibre sequence
 \begin{equation*}
 {\bf R}\Gamma_{\EU N}(k^\circ,M^\circ) \lto  
 {\bf R}\Gamma_{\EU N}(k,M)\lto 
 {\bf R}\Gamma_{\EU N}(k^\circ,M^{\circ,\lozenge}(1)[1])^* \overset{+1}{\lto}
 \end{equation*}
\end{proposition}

\begin{proof} 
We immediately reduce to the case where $M^\circ$ is static. Then
\begin{align*}
{\bf R}\Gamma_{\EU N}(k^\circ,M^\circ) &= {\rm fib}\left[M^{\circ,0} \overset{\nabla}{\lto}
\widehat{\Omega}^1_{C^\circ} \otimes_{C^\circ} M^{\circ,1} \right] \\
{\bf R}\Gamma_{\EU N}(k,M) &= {\rm fib}\left[M^{0} \overset{\nabla}{\lto}
\widehat{\Omega}^1_{C} \otimes_{C} M^1 \right]
\end{align*}
and so by Lemma \ref{lemma: pont dual orthog} the cofibre of ${\bf R}\Gamma_{\EU N}(k^\circ,M^\circ)\to {\bf R}\Gamma_{\EU N}(k,M)$ is given by the totalisation of the complex
\[
\left(\widehat{\Omega}^1_{C^\circ} \otimes_{C^\circ}  M^{\circ,0,\vee}\right)^* \lto \left(  M^{\circ,1,\vee}\right)^*
\]
with terms sitting in horizontal cohomological degrees $0$ and $1$. But we can identify this with
\[
 \left( M^{\circ,\lozenge,-1} \lto  \widehat{\Omega}^1_{C^\circ} \otimes_{C^\circ}  M^{\circ,\lozenge,0}\right)^*
\]
where the complex inside the brackets sits in horizontal cohomological degrees $-1$ and $0$, and is hence nothing other than ${\bf R}\Gamma_{\EU N}(k^\circ,M^{\circ,\lozenge}(1)[1])$.
\end{proof}
 
 We can now complete the proof of Theorem \ref{theo: lagrange}. Take the fibre sequence
 \[
 M^{\circ,\vee}(1) \lto
 j^{(0)}_+(M^{\circ,+\infty,\vee})\oplus j^{(-1)}_+(M^{\circ,-\infty,\vee}) \lto
 M^{\circ,\lozenge}(1)\overset{+1}{\lto}
 \]
 and apply Nygaard cohomology to give a fibre sequence
 \[
 {\bf R}\Gamma_{\EU N}(k^\circ,M^{\circ,\vee}(1)) \lto
 {\bf R}\Gamma_{{\rm dR},(0)}(k^\circ,M^{\circ,\vee}) \oplus {\bf R}\Gamma_{{\rm dR},(-1)}(k^\circ,M^{\circ,\vee}) \lto
 {\bf R}\Gamma_{\EU N}(k^\circ,M^{\circ,\lozenge}(1))\overset{+1}{\lto}
 \]
Exactly as in the proof of Theorem \ref{theo: main1}, we can pull this back along the `diagonal' map from ${\bf R}\Gamma_{{\rm dR},(0)}(k^\circ,M^{\circ,\vee})$ to the central term to give a fibre sequence
 \[
 {\bf R}\Gamma_{\rm Syn}(k^\circ,M^{\circ,\vee}(1)) \lto
 {\bf R}\Gamma_{{\rm dR},(0)}(k^\circ,M^{\circ,\vee}) \lto
 {\bf R}\Gamma_{\EU N}(k^\circ,M^{\circ,\lozenge}(1))\overset{+1}{\lto}
 \]
 which after shifting, dualising, and rotating gives a fibre sequence
 \[ {\bf R}\Gamma_{\rm Syn}(k^\circ,M^{\circ,\vee}(1)[2])^* \lto  {\bf R}\Gamma_{\EU N}(k^\circ,M^{\circ,\lozenge}(1)[1])^* \lto {\bf R}\Gamma_{{\rm dR},(0)}(k^\circ,M^{\circ,\vee}[1])^*  \overset{+1}{\lto} \]
We now use the standard exact sequence computing syntomic cohomology, together with Proposition \ref{prop: N lagrange}, to form the commutative diagram
\[
\xymatrix{
{\bf R}\Gamma_{\rm Syn}(k^\circ,M^\circ)\ar[r]\ar[d] & {\bf R}\Gamma_{\EU N}(k^\circ,M^\circ) \ar[r]\ar[d] & {\bf R}\Gamma_{{\rm dR},(0)}(k^\circ,M^\circ)\ar[r]^-{+1}\ar[d] & \\
 {\bf R}\Gamma_{\rm Syn}(k,M)\ar[r] \ar[d] & {\bf R}\Gamma_{\EU N}(k,M) \ar[r]\ar[d] & {\bf R}\Gamma_{{\rm dR},(0)}(k,M)\ar[r]^-{+1}\ar[d] & \\ 
 {\bf R}\Gamma_{\rm Syn}(k^\circ,M^{\circ,\vee}(1)[2])^* \ar[d]^{+1}\ar[r] &  {\bf R}\Gamma_{\EU N}(k^\circ,M^{\circ,\lozenge}(1)[1])^* \ar[r]\ar[d]^{+1} & {\bf R}\Gamma_{{\rm dR},(0)}(k^\circ,M^{\circ,\vee}[1])^* \ar[d]^{+1}\ar[r]^-{+1} & \\
   & & & 
   }\]
in which all rows, as well as the middle and right hand columns, are fibre sequences. It follows that the left hand column is also a fibre sequence, which finishes the proof.

\begin{remark} If $p>2$, then \cite[Theorem 1.1]{TV24} also applies over $k^\circ$, and as explained in \cite[p.662, final paragraph]{BKK26} the restriction $p>2$ is unnecessary. Hence the fact that the analogue of Theorem \ref{theo: syn cohom complex 2} holds over $k^\circ$ means that the analogue of Corollary \ref{cor: syn coh flat} also holds over $k^\circ$. In particular, we can recover the classical orthogonality results for flat cohomology (that is, \cite[Corollary III.7.2]{Mil86}) from Theorem \ref{theo: lagrange}, at least once we know injectivity of the map
\begin{equation} \label{eqn: int flat coh} \rH^*_{\rm Syn}(k^\circ,\fr{g}^{(0)}H)\to \rH^*_{\rm Syn}(k,\fr{g}^{(0)}H)\end{equation}
whenever $H$ is a finite flat group scheme (of $p$-power order) over $k^\circ$. But, given the equivalence with flat cohomology, the argument that \eqref{eqn: int flat coh} is an isomorphism in degree $0$ and injective in degree $1$ is completely elementary, and this is then enough to show vanishing of the left hand side of \eqref{eqn: int flat coh} in all higher degrees by self-duality of the fibre sequence \eqref{eqn: fibre}. Thus we obtain the required injectivity.
\end{remark}

\bibliographystyle{bibsty}
\bibliography{../../Templates/lib.bib}

\end{document}

%% file: headart.tex
\usepackage[toc,page]{appendix}
\usepackage[headheight=12.1pt]{geometry}
\usepackage[all,cmtip]{xy}
\usepackage{amsmath}
\usepackage{amsfonts}
\usepackage{amssymb}
\usepackage{mathrsfs}   
\usepackage{amsthm}

\usepackage{mathtools}
\usepackage{xcolor}
\usepackage{stmaryrd}
\usepackage{fancyhdr}
\usepackage[backref=page]{hyperref}

\usepackage{euscript}
\newcommand{\EU}[1]{\EuScript{#1}}

\usepackage{amsbsy}
\SetSymbolFont{stmry}{bold}{U}{stmry}{m}{n}
\usepackage{bm}

\hypersetup{
     colorlinks   = true,
     citecolor    = red,
     linkcolor = blue,
     urlcolor = purple
}
    
\numberwithin{equation}{section}
   
\makeatletter
\newcommand*{\relrelbarsep}{.386ex}
\newcommand*{\relrelbar}{%
  \mathrel{%
    \mathpalette\@relrelbar\relrelbarsep
  }%
}
\newcommand*{\@relrelbar}[2]{%
  \raise#2\hbox to 0pt{$\m@th#1\relbar$\hss}%
  \lower#2\hbox{$\m@th#1\relbar$}%
}
\providecommand*{\rightrightarrowsfill@}{%
  \arrowfill@\relrelbar\relrelbar\rightrightarrows
}
\providecommand*{\leftleftarrowsfill@}{%
  \arrowfill@\leftleftarrows\relrelbar\relrelbar
}
\providecommand*{\xrightrightarrows}[2][]{%
  \ext@arrow 0359\rightrightarrowsfill@{#1}{#2}%
}
\providecommand*{\xleftleftarrows}[2][]{%
  \ext@arrow 3095\leftleftarrowsfill@{#1}{#2}%
}
\makeatother

\theoremstyle{plain}
\newtheorem*{theoremu}{Theorem}
\newtheorem{theorem}{Theorem}[section]

\newtheorem{proposition}[theorem]{Proposition}

\newtheorem{corollary}[theorem]{Corollary}

\newtheorem{lemma}[theorem]{Lemma}

\theoremstyle{definition}

\newtheorem{definition}[theorem]{Definition}

\theoremstyle{remark}
\newtheorem{remark}[theorem]{Remark}

\newtheorem{example}[theorem]{Example}

\newcommand{\N}{{\mathbb N}}
\newcommand{\Z}{{\mathbb Z}}
\newcommand{\Q}{{\mathbb Q}}

\newcommand{\F}{{\mathbb F}}

\newcommand{\G}{{\mathbb G}}

\newcommand{\ca}[1]{\mathcal{#1}}
\newcommand{\fr}[1]{\mathfrak{#1}}
\newcommand{\rH}{{\rm H}}
\newcommand{\bD}{{\bf D}}

\newcommand{\colim}[1]{\underset{#1}{\mathrm{colim}}\,}
\renewcommand{\lim}[1]{\underset{#1}{\mathrm{lim}}\,}
\newcommand{\isomto}{\overset{\cong}{\longrightarrow}}

\newcommand{\from}{\colon}
\renewcommand{\to}{\rightarrow}

\newcommand{\lto}{\longrightarrow}

\newcommand{\lser}[1]{(\!(#1)\!)}
\newcommand{\pow}[1]{\llbracket #1 \rrbracket}

\newcommand{\spf}[1]{\mathrm{Spf}\!\left(#1\right)}
\newcommand{\dR}{\mathrm{dR}}
\newcommand{\Tr}{\mathrm{Tr}}
\newcommand{\Res}{\mathrm{Res}}

\renewcommand{\d}{{\rm d}}
\newcommand{\dlog}{{\rm d}\!\log}